\documentclass[12pt]{article}
\usepackage{amsmath,amsthm,verbatim,amssymb,amsfonts,amscd,diagrams, graphics}
\theoremstyle{plain}
\newtheorem{theorem}{Theorem}[section]
\newtheorem{corollary}[theorem]{Corollary}
\newtheorem{lemma}[theorem]{Lemma}
\newtheorem{proposition}[theorem]{Proposition}

\theoremstyle{definition}
\newtheorem{definition}[theorem]{Definition}

\theoremstyle{remark}
\newtheorem{remark}[theorem]{Remark}

\newarrow{ul}---->
\newarrow{Backwards}<----   

\newcommand{\onto}{\twoheadrightarrow}
 
\newcommand{\Td}[1]{\Tilde{#1}}
\newcommand{\td}[1]{\tilde{#1}}
\newcommand{\into}{\hookrightarrow}
\newcommand{\Z}{\mathbb{Z}}
\newcommand{\Q}{\mathbb{Q}}
\newcommand{\R}{\mathbb{R}}
\newcommand{\bd}{\partial}

\newcommand{\mc}[1]{\mathcal{#1}}

\newcommand{\Hom}{\text{Hom}}
\newcommand{\Ext}{\text{Ext}}
\newcommand{\mf}{\mathfrak}

\newcommand{\im}{\text{im}}
\newcommand{\cok}{\text{cok}}
\newcommand{\coim}{\text{coim}}

\begin{document}

\title{Cobordism of disk knots}
\author{Greg Friedman\\Yale University}
\date{January 14, 2004}
\maketitle
typeset=\today

\begin{abstract}
We study cobordisms and cobordisms rel boundary of PL locally-flat disk knots $D^{n-2}\into D^n$. Cobordisms of disk knots that do not fix the boundary sphere knots are easily classified by the cobordism properties of these boundaries, and any two even-dimensional disk knots with isotopic boundary knots are cobordant rel boundary. However, the cobordism rel boundary theory of odd-dimensional disk knots is more subtle. Generalizing results of Levine on cobordism of sphere knots, we define disk knot Seifert matrices and show that two higher-dimensional disk knots with isotopic boundaries are cobordant rel boundary if and only if their disk knot Seifert matrices are algebraically cobordant. We also find necessary and sufficient conditions to realize a  Seifert matrix cobordism class among the disk knots corresponding to a fixed boundary knot, assuming the boundary knot has no middle-dimensional $2$-torsion. This classification is performed by relating the Seifert matrix of a disk knot to its Blanchfield pairing and by establishing a close connection between this Blanchfield pairing and the Farber-Levine torsion pairing of the boundary knot (in fact, for disk knots satisfying certain connectivity assumptions, the disk knot Blanchfield pairing will determine the boundary Farber-Levine pairing). In asddition, we  study the dependence of disk knot Seifert matrices on choices of Seifert surface, demonstrating that all such Seifert matrices are \emph{rationally} S- equivalent, but not necessarily integrally S-equivalent.

\end{abstract}


\tableofcontents

\section{Introduction}

Two locally-flat sphere knots, i.e. PL locally-flat embeddings $K_0, K_1: S^{n-2}\into S^n$, are called cobordant (sometimes concordant) if there exists a proper locally-flat PL embedding $\mf{K}: S^{n-2}\times [0,1]\into S^n\times [0,1]$ such that $\mf{K}|S^{n-2}\times 0=K_0$ and  $\mf{K}|S^{n-2}\times 1=-K_1$. Here $-K_1$ is the knot obtained by $K_1$ followed by a reflection of $S^n$. This ``negative'' knot occurs due to the usual  reversal of orientation at one end of a cobordism. It was shown by Kervaire \cite{Ke} that all knots of even dimension are cobordant to the trivial unknotted embedding, so all even-dimensional knots are \emph{null-cobordant}, or \emph{slice}. For odd-dimensional knots of  dimension  $n>3$, Levine \cite{L69} obtained complete necessary and sufficient algebraic conditions for two knots  to be cobordant. These conditions are stated in terms of the Seifert matrices of knots, and two knots are cobordant if and only if their Seifert matrices satisfy a relationship of algebraic cobordism (a knot does not determine a unique Seifert matrix, but any two such Seifert matrices will lie in the same algebraic cobordism class). Levine also demonstrates the existence of knots that realize any possible algebraic cobordism class, within the other restrictions necessary for a matrix to be a Seifert matrix. 

In this paper, we turn our attention to the cobordism of disk knots, PL locally-flat proper embeddings $D^{n-2}\into D^n$. Since the embeddings are proper, each disk knot $L$ determines a locally-flat sphere knot $K$ on restriction to the boundary. We will call two disk knots $L_0, L_1$ cobordant  if there exists a proper locally-flat PL embedding $\mf{L}: D^{n-2}\times [0,1]\into D^n\times [0,1]$ such that $\mf{L}|D^{n-2}\times 0=L_0$ and  $\mf{L}|S^{n-2}\times 1=-L_1$. Note that the restriction of $\mf{L}$ to $\bd D^{n-2}\times [0,1]$ provides a cobordism between the boundary sphere knots $K_0$ and $K_1$. If this cobordism extends to an ambient isotopy of  $K_0$ to $K_1$, we will call $\mf{L}$ a cobordism rel boundary. 

The cobordism theory of sphere knots was first studied by  Fox and Milnor  \cite{M66} as they sought to remove singularities of embeddings of manifolds by replacing cones on smooth knots by ``slicing disks''. By studying cobordism rel boundary of disk knots, we seek to classify  precisely  such slicing disks up to their own cobordisms, so in some sense we are studying a second order of cobordism theory. The results of this theory will provide some measure of the number of ways in which a codimension  two  embedding with point singularities of a manifold  can be converted into a smooth embedding via the local resolution of singularities. On the other hand, there is another close relation between smooth disk knots and sphere knots with point singularities (see \cite{GBF1}), and the cobordism theory of the former will be essential in studying that of the latter. In a future paper, we will study the cobordisms of knots with point singularities and obtain a similar measure of the number of ways to remove $0$-dimensional strata from manifold embeddings with $1$-dimensional singularities (the natural goal is to  reduce the number of necessarily distinct strata in a stratification). 

We now outline the main results of this paper. Three of the four cases of interest can be studied rather easily, and we will obtain the following results almost immediately:

\begin{proposition}[Proposition \ref{P: even cob rel}]
If $n$ is even, then two disk knots $L_0, L_1: D^{n-2}\to D^n$ are cobordant if and only if their boundary knots are cobordant. 
\end{proposition}

\begin{proposition}[Proposition \ref{P: even cob}]
If $n$ is even, then any two disk knots $L_0, L_1: D^{n-2}\to D^n$ with isotopic boundary knots are  cobordant rel boundary. 
\end{proposition}

\begin{proposition}[Proposition \ref{P: odd cob}]
If $n$ is odd, then any two knots $L_0, L_1: D^{n-2}\to D^n$  are cobordant.
\end{proposition}

This leaves the more challenging case of cobordism rel boundary for odd dimensional disk knots. To study this case, we will need to introduce Seifert matrices for disk knots. As opposed to Seifert matrices for sphere knots, which arise as certain forms  on the middle dimensional homology of Seifert surfaces, Seifert matrices for disk knots are forms defined only on certain quotient homology modules. Disk knot Seifert matrices also may differ from those for sphere knots in that, if $A$ is such a Seifert matrix, the matrix $A+(-1)^nA'$ need not be integrally unimodular, only rationally so. Nonetheless, algebraic cobordism is well-defined on this larger class of matrices, and we attain the following conclusion:

\begin{theorem}[Theorem \ref{T: main cobordism}]
Let $L_0, L_1: D^{2n-1}\into D^{2n+1}$, $n>1$, be two disk knots with the same boundary knot.
Let $A_0$ and $A_1$ be Seifert matrices for  
 $L_0$ and $L_1$, respectively. Then $L_0$ and $L_1$  are cobordant rel boundary if and only if $A_0$ and $ A_1$ are cobordant.
\end{theorem}

Several interesting corollaries  follow:

\begin{corollary}[Corollary \ref{C: Alex}]
Suppose that $L_0$ and $L_1$ are disk knots $D^{2n-1}\subset D^{2n+1}$, $n>1$, such that $\bd L_0=\bd L_1=K$. Then a necessary condition for $L_0$ and $L_1$ to be cobordant rel boundary is that the product of the middle-dimensional Alexander polynomials $c^{L_0}_{n}(t)c_{n}^{L_1}(t)$ be similar in $\Q[t,t^{-1}]$ to a polynomial of the form $p(t)p(t^{-1})$.
\end{corollary}

\begin{theorem}[Theorem \ref{T: there is a sphere}]
Let $L_0$ and $L_1$ be two disk knots $D^{2n-1}\subset D^{2n+1}$, $n>1$, with common boundary $K$. Then there exists a sphere knot $\mc{K}: S^{2n-1}\subset S^{2n+1}$ such that $L_0$ is cobordant to the knot sum (away from the boundary) $L_1\# \mc {K}$. 
\end{theorem}

\begin{theorem}[Theorem \ref{T: cob to simple}]
Given any disk knot $L:D^{2n-1}\subset D^{2n+1}$, $n>1$, $L$ is cobordant rel boundary to a disk knot $L_1$ such that $\pi_i(D^{2n+1}-L_1)\cong \pi_i(S^1)$ for $i<n$.
\end{theorem}

This last theorem tells us that every disk knot is cobordant rel boundary to a \emph{simple disk knot}.

The next question to consider is that of which cobordism classes of matrices arise as the Seifert matrices of disk knots. We will show that all possible such matrices occur for knots of sufficiently high dimension, but we will also be interested in the sharper question of which classes arise for disk knots given a fixed boundary knot. It turns out that if the boundary knots has no middle-dimensional $2$ -torsion then the determining information from the boundary knot is its Farber-Levine torsion pairing since, as we will see, there is a close relationship between the  Blanchfield pairing of a disk knot, which is determined by its Seifert matrices, and the Farber-Levine torsion pairing of its boundary knot. In fact, if a disk knot is simple, i.e. its complement has the homotopy groups of a circle below the ``middle'' dimension, its Blanchfield pairing will completely determine the Farber-Levine pairing of its boundary knot. In particular, $\td C$ is the infinite cyclic cover of the disk knot complement and $\td X$ is the infinite cyclic cover of the complement of its boundary sphere knot, we can prove the following:

\begin{theorem}[Theorem \ref{T: Blanch gives FL}]
Given a simple disk knot $D^{2n-1}\into D^{2n+1}$, the module $H_{n-1}(\td X)$ and the Farber-Levine $\Z$-torsion pairing on its $\Z$-torsion submodule  $T_{n-1}(\td X)$  are determined up to isometry by the isometry class of the Blanchfield self-pairing on $H_n(\td C)$. 
\end{theorem}

\begin{theorem}[Theorem \ref{H to FL}]
For a simple disk knot $L:D^{2n-1}\subset D^{2n+1}$, the $\Lambda$-module $T_{n-1}(\td X)$ and its Farber-Levine $\Z$-torsion pairing are determined up to isometry by the isometry class of $\cok(H_n(\td X)\to H_n(\td C))$ with its self-Blanchfield pairing.
\end{theorem}

\begin{corollary}[Corollary \ref{C: Seifert gives FL}]
For a simple disk knot $L:D^{2n-1}\subset D^{2n+1}$, the $\Lambda$-module $T_{n-1}(\td X)$ and its Farber-Levine $\Z$-torsion pairing are determined up to isometry by any Seifert matrix for $L$.
\end{corollary} 

These theorems, together with  a theorem of Kojima \cite{Koj}, will allow us to prove that, given a fixed boundary knot  $K$ of sufficiently high dimension and with no middle-dimensional $2$-torsion, any cobordism class of matrices containing an element that correctly determines the Farber-Levine pairing of $K$ is realizable as cobordism class of  Seifert matrices of a  disk knot with $K$ as its boundary knot. See Theorem \ref{T: realize} for a more accurate statement.

In the course of these investigations, we will also need to engage in an in-depth study of  how the Seifert matrix of a disk knot varies with choice of Seifert surface. In particular, in an extended technical section we will prove the following theorem:

\begin{theorem}[Theorem \ref{T: s-equiv}]
Any two Seifert matrices for a disk knot differ by a rational S-equivalence. 
\end{theorem}

The organization of this paper is as follows: In Section \ref{S: basics}, we present the basic technical details concerning Seifert matrices of disk knots. In Section \ref{S: cobordism}, we begin our geometric investigation and determine when two disk knots are cobordant. Section \ref{S: realize} contains the constructions that allow us to realize the algebraic cobordism matrices geometrically.  Sections \ref{S: blanch to F-L} contains the discussion of the relation between disk knot Blanchfield pairings and their boundary sphere knot Farber-Levine pairings. Finally, Section \ref{S: seifert surfaces} contains the calculations of how disk knot Seifert matrices change as the Seifert surface is varied.

\section{Seifert matrix basics}\label{S: basics}

We first introduce some notation that will be in constant use: Let $L$ denote a PL locally-flat  \emph{disk knot} $L: D^{n-2}\into D^n$. All disk knots will be proper embeddings, i.e. $D^{n-2}\cap \bd D^n=\bd D^{n-2}$, and there is a collar of the boundary in which the embedding is PL-homeomorphic to $(\bd D^n, \bd D^{n-2})\times I $. The boundary embedding $\bd D^{n-2}\into \bd D^n$ is the locally-flat boundary sphere knot $K$.  We will employ the standard abuse of notation and confuse the maps $L$ and $K$ with their images. We use $C$ to denote the \emph{exterior} of $L$, the complement of an open regular neighborhood of $L$; $C$ is homotopy equivalent to $D^n-L$. We use $X$ to denote $C\cap \bd D^n$, the exterior of $K$. Using Alexander duality (respectively, Alexander duality for a ball (see \cite[p. 426]{MK}), $X$ and $C$ are homology circles and so possess infinite cyclic covers that we denote $\td X$ and $\td C$. $F$ denotes a Seifert surface for $K$, and $V$ denotes a \emph{Seifert surface} for $L$, i.e. an oriented bi-collared $n-1$-dimensional submanifold of $D^n$ whose boundary is the union of $L$ and a Seifert surface for the boundary knot $K$. Such Seifert surfaces always exist (see \cite{GBF1}). Note that $H_*(\bd V)\cong H_*(F)$ for $*\leq n-3$. 

The groups $H_*(\td X)$, $H_*(\td C)$, and $H_*(\td C,\td X)$ inherit structures as modules over $\Lambda=\Z[\Z]=\Z[t,t^{-1}]$ by the action of the covering translation. A $\Lambda$-module is of \emph{type K} if it is finitely generated and multiplication by $t-1$ acts as an automorphism. Equivalently, a  $\Lambda$-module of type K is a finitely generated $\Lambda[(t-1)^{-1}]$ module. It is well known that $H_*(\td X)$ is a torsion $\Lambda$-module of type K for $*>0$ (see e.g. \cite{L77}).  Since $C$ is a homology circle, $H_*(\td C)$ is  also of type K for $*>0$ by Levine \cite[Prop. 1.2]{L77} since the proof of this proposition only relies on $C$ being a homology circle. It then follows from \cite[Cor. 1.3]{L77} that $H_*(\td C)$, $*>0$, is a $\Lambda$-torsion module. Hence so is $H_*(\td C, \td X)$ from the reduced long exact sequence of the pair (in fact, it is similarly of type K by the five lemma applied to the long exact sequence of the pair under multiplication by $t-1$).

Let $V$ be a Seifert surface of a knot $L:D^{2n-1}\into D^{2n+1}$ with boundary Seifert surface $F$. Then $H_i(V,F)\cong H_i(V, \bd V)$ for $i\leq 2n-2$ and $H_{2n-1}(V,F)\to H_{2n-1}(V, \bd V)$ is onto. So, in particular, $H_n(V,F)\cong H_n(V, \bd V)$, induced by inclusion, for $n\geq 2$. Poincar\'e-Lefschetz duality implies that we have a nonsingular intersection pairing $S: F_n(V)\otimes F_n(V, F)\to \Z$, where $F_i$ is $H_i$ modulo torsion. In particular, the ranks of $F_n(V)$ and $F_n(V, F)$ are equal.

We next need to investigate the duality properties of these modules more carefully around the middle dimension. 
Consider the portion of the long exact sequence of the pair given by
\begin{equation*}
\begin{CD}
H_n(F)&@>i_*>>& H_n(V)&@>p_*>>& H_n(V,F)&@>\bd_*>>& H_{n-1}(F).
\end{CD}
\end{equation*}
Let $E$ be the kernel of $\bd_*$ mod torsion, and let $\bar E$ be the cokernel of $i_*$ mod torsion. If $E\cong Z^m$ and $F_n(V, F)\cong \Z^k$, then it is possible to choose a basis of $F_n(V,F)$ so that $E$ is a subgroup of the subgroup $\td E\subset \Z^k$ consisting of the first $m$ $\Z$ summands of $F_n(V,F)$.  This follows from the existence of a diagonal matrix representing $p_*$ mod torsion (see \cite[Thm. 11.3]{MK}. In fact, we can further assume by this method that there are  generators $\alpha_i$ of $\Z$ summands of $\td E$ and non-zero least integers $q_i>0$ such that $q_i\alpha_i$ is in the image of $p_*$ (mod torsion). In other words, we can assume that $$p_*=\begin{pmatrix} q_1 \\
 &q_2\\
&&\ddots \\
&&&q_m\\
&&&&0\\
&&&&&\ddots
\end{pmatrix}.$$ 
Here $\td E$ is the subgroup represented by the first $m$ columns of this matrix.

 In what follows, we will use $i_*$, $p_*$, and $\bd_*$ also to denote the maps induced mod torsion. 

To fix notation, let $\{\alpha_i\}_{i=1}^k$ be a basis for $F_n(V,F)$ and $\{\alpha_i\}_{i=1}^m$ the restricted basis for $\td E$. Let $\{\delta_i\}_{i=1}^k$ be the dual basis in $F_n(V)$ under $S$, i.e. $S(\alpha_i,\delta_j)$ equals $1$ if $i=j$ and $0$ if $i\neq j$. We want to show that $\{\bar \delta_i\}_{i=1}^m$ is a basis for $\bar E$, where $\bar \delta_i$ is the projection of $\delta_i$ to $\bar E$. 

Let $\bar S: F_n(V)\otimes F_n(V)\to \Z$ be the pairing defined by $\bar S(a,b)=S(p_*(a), b)$. Note that this is simply equal to the intersection number of chains representing $a$ and $b$, and so also $\bar S(a,b)=(-1)^nS( p_*(b),a )$. This induces a pairing $T: \bar E\otimes \bar E\to \Z$. To see that this is well-defined, we need to show that $\bar S(a,b)=0$ if $a$ or $b$ lies in $\ker p_*$, but this is evident. 

We now claim that the dual elements to $\{\alpha_i\}_{i=1}^m$ under $S$ is a basis of $\bar E$ under the canonical projection of $F_n(V)$ onto $\bar E$:

\begin{proposition}
The elements $\{\bar \delta_i\}_{i=1}^m$ are a basis for $\bar E$.
\end{proposition}
\begin{proof}

We must have $\delta_i\notin \ker p_*$ for $i\leq m$: If $\delta_i\in \ker p_*$, then for any $\eta\in F_n(V)$, $\bar S(\delta_i, \eta)=0$, as noted above. But we know that $S(\delta_i, \alpha_i)=1$ and that, by our choice of the basis $\{\alpha_i\}$, there exists an integer $q_i\neq 0$ and a $\chi_i\in F_n(V)$ such that $q_i\alpha_i= p_*(\chi_i)$. So if $\delta_i\in \ker(p_*)$, then $q_i=S(\delta_i, q_i\alpha_i)= \bar S(\delta_i, \chi_i)=0$, a contradiction. Thus $\delta_i\notin \ker p_*$. 
Furthermore, there is no non-trivial linear combination $v=\sum d_i\delta_i\in \ker p_*$, else similarly $q_id_i=S(v, q_i\alpha_i)=\bar S(v, \chi_i)=0$. Thus the $\bar \delta_i$ are linearly independent in $\bar E$.

Now, $\bar E$ and $\td E$ must have the same rank. In fact, $E$ and $\bar E$  are isomorphic, being the respective image and coimage of $p_*$, and $E$ and $\td E$ have the same rank $m$ by the construction of $\td E$. It remains to show that the $\bar \delta_i$, $i\leq m$, span integrally.

Let $\bar x$ be an element of $\bar E$, let $\pi: F_n(V)\to \bar E\cong F_n(V)/(\ker (p_*)$ be the projection (here $\ker(p_*)$  is assumed to have had its torsion quotiented out already), and let $x\in F_n(V)$ such that $\pi(x)=\bar x$. Then $x=\sum_{i=1}^kn_i\delta_i$, $n_i\in \Z$. Similarly, let $\pi_{\Q}:F_n(V)\otimes \Q \to (F_n(V)\otimes \Q)/(\ker(p_*)\otimes \Q)=\bar E\otimes \Q$ be the  projection and let $y\in F_n(V)\otimes \Q$ such that $\pi_{\Q}(y)=\bar x\otimes 1$. Since we know that $\bar E\otimes \Q$ is spanned by $\{\bar \delta_i\otimes 1\}_{i=1}^m$, we can choose a  $y$ of the form $y=\sum_{i=1}^m r_i(\delta_i\otimes 1)$, $r_i\in \Q$. Since $\pi(x)\otimes 1=\pi_{\Q}(x\otimes 1)=\pi_{\Q}(y)$, $x\otimes 1-y\in \ker(\pi_{\Q})$. Now let $S_{\Q}$ and $\bar S_{\Q}$ denote the rational pairings induced from $S$ and $\bar S$. Then for $1\leq j\leq m$, 
\begin{align*}
0&=S_{\Q}(q_j\alpha_j \otimes 1,x\otimes 1-y))\\
&=S_{\Q}(p_*(\chi_j\otimes 1),x\otimes 1-y)\\
&=\bar S_{\Q}(\chi_j\otimes 1, x\otimes 1-y)\\
&=\pm \bar S_{\Q}(x\otimes 1-y, \chi_j\otimes 1)\\
&=\pm S_{\Q}(p_*( x\otimes 1-y), \chi_j\otimes 1).
\end{align*}
But $S_{\Q}(q_j\alpha_j\otimes 1 ,x\otimes 1-y)=q_jS_{\Q}(\alpha_j \otimes 1, \sum_{i=1}^kn_i(\delta_i\otimes 1)-\sum_{i=1}^m r_i(\delta_i\otimes 1))=q_j(n_j-r_j)$ since the $\alpha_i$ and $\delta_i$ remain dual bases rationally. Thus we see that $n_j=r_j$. Therefore,  $x\otimes 1-y=\sum_{i=m+1}^k n_i\delta_i$. But this is in $\ker (\pi_{\Q})\cap F_n(V)=\ker (\pi)$. So $\bar x=\pi(x)=\pi(x- \sum_{i=m+1}^k n_i\delta_i)=\pi(\sum_{i=1}^m n_i\delta_i)$. This shows that  $\bar x$ is in the integral span $\{\bar \delta_i\}_{i=1}^m$ and completes the proof.
\end{proof}

\begin{corollary}
For $i>m$, $\delta_i\in \ker (p_*)$.
\end{corollary}
\begin{proof}
It follows from the proof of the preceding proposition by taking $x=\delta_i$ that $\pi_*(\delta_i)=0$ for $i>m$ and hence $\delta_i\in \ker (p_*)$.
\end{proof}

Now consider again $p_*:F_n(V)\to F_n(V,F)$. If $\bar \delta_i$, $i\leq m$ is a basis element of $\bar E$ as above, then $p_*(\bar \delta_i)=\sum_{i=1}^m R_{li}\alpha_l$, and $R_{ji}=S( p_*(\bar \delta_i),\delta_j)=\bar S(\delta_i,  \delta_j)=T(\bar \delta_i, \bar \delta_j)$, since we have noted that the intersection pairing is trivial on elements in the kernel of $p_*$. Thus with these bases, the transpose of the  matrix $R$ of the mapping $p_*$ is the matrix of the intersection pairing $T$ on $\bar E$. In other words, we have proven the following:

\begin{proposition}
The matrix $R_{ji}$ of the mapping $p_*:\bar E\to \td E$ with respect to dual bases is the matrix of the intersection pairing $T(\bar \delta_i,\bar \delta_j$.

\end{proposition}

At this point, we note that since $\td E$ is a well-defined subspace of $F_n(V,F)$, the $\{\delta_i\}$, $i\leq m$, span a well-defined dual subspace in $F_n(V)$. We have already noted that $\bar E$ shares this basis set under the projection from $F_n(V)$. Hence to simplify notation below, we will identify $\bar E$ with the subspace of $F_n(V)$ spanned by the $\{\delta_i\}$, $i\leq m$, and remove the bars from the $\delta$ notation. We can also then consider $T$ as the restriction of $S$ to this subspace. 

\vskip.5cm

Similarly consideration to those above occur on the complement of $V$ in $D^{2n+1}$. We let $Y=D^{2n+1}-V$ and $Z=\bd D^{2n+1}-F$. Then we can use the map $p^Y_*:F_n(Y)\to F_n(Y,Z)$ to define $G$, $\td G$, and $\bar G$ analogously to  $E$, $\td E$, and $\bar E$. It follows from Alexander duality for a ball and the computations in \cite[\S 3.6.3]{GBF1} that the ranks of $G$, $\td G$, and $\bar{\Td G}$ will also be $m$.

Also by the arguments in \cite[\S 3.6.3]{GBF1}, which are  similar to those above, we can now also take as a basis of $\bar G$ the duals $\{\beta_i\}_{i=1}^m$ such that $L'(\alpha_i, \beta_j)=\delta_{ij}$, where $L': F_n(V, F)\otimes F_n(Y)\to \Z$ is the Alexander linking pairing for a ball (see \cite{GBF1}). For $\td G$, assume that a basis $\{\gamma'_i\}_{i=1}^k$ is chosen so that $p_*^Y$ can be diagonalized with $\td G$ in this basis and such that $G$ lies in the span of $\{\gamma'_i\}_{i=1}^m$. We know also from  \cite{GBF1} that if the duals $\{\delta'_i\}_{i=1}^k$ such that $L''(\gamma'_i, \delta'_j)= \delta_{i,j}$ are another basis for $F_n(V)$ (where $L'': F_n(Y, Z)\otimes F_n(V)\to \Z$ is the other associated linking pairing), then the projections $\{\bar \delta_i'\}_{i=1}^m$ of  $\{\delta'_i\}_{i=1}^m$ are also a basis for $\bar E$ and $\{\delta'_i\}_{i=m+1}^k$ is a basis for $\ker p_*$.  In particular then, we can change the basis $\{\gamma'_i\}_{i=1}^m$ to a basis $\{\gamma_i\}_{i=1}^m$ that is dual to $\{\delta_i\}_{i=1}^m$: if $\Xi$ is a change of basis matrix such that $\bar \delta_i=\sum_{j=1}^m\Xi_{ij}\bar \delta'_j$ and $\Theta=(\Xi^{-1})'$ (here $'$ indicates transpose), then let $\gamma_i=\sum_j\Theta_{ij}\gamma'_j$. So then $L''(\gamma_i, \delta_j)=L''( \sum_a\Theta_{ia}\gamma'_a, \sum_{b}\Xi_{jb}\delta'_b)=\sum_a \Theta_{ia}L''(\gamma'_a, \sum_{b}\Xi_{jb}\delta'_b)=\sum_a \Theta_{ia} \Xi_{ja}=\delta_{i,j}$. The first equality holds since $L''(x,y)=0$ if $y\in \ker(p_*)$ and $x\in \td G$ (see \cite[\S 3.6.3]{GBF1}). 

In other words, we have simply proven the following lemma:

\begin{lemma} We can choose bases for $\bar E$, $\bar G$, $\td E$, and $\td G$ so that those for $\bar E$ and $\td E$ are dual under the intersection pairing $S$ and those for $\bar E$ and $\td G$ and those for $\td E$ and $\bar G$ are dual under the linking pairings $L''$ and $L'$.  Furthermore, the basis for $\td E$ is one with respect to which the map $p_*$ can be diagonalized by changing the basis of $\bar E$.
\end{lemma}

Now, as in \cite{GBF1},  we let $i_{+*}$ and $i_{-*}$ denote the maps on homology induced by pushing $V$ off itself along the bicollar in the respective negative and positive directions (note the reversed order; we follow the convention of \cite{L66}). With respect to the above bases, we can define matrices $\lambda$, $\sigma$, $\tau$, and $\mu$ by
\begin{align*}
i_{+*}(\delta_j)&=\sum_i \lambda_{ij} \beta_i\\
i_{-*}(\delta_j)&=\sum_i \sigma_{ij} \beta_i\\
i_{+*}(\alpha_j)&=\sum_i \mu_{ij} \gamma_i\\
i_{-*}(\alpha_j)&=\sum_i \tau_{ij} \gamma_i,
\end{align*}
where all indices run from $1$ to $m$. It is shown in \cite{GBF1} that the  matrix $L''(i_{-}p_*(\delta_i), \delta_j)=\sum_{k=1}^m R_{ki}\tau_{jk}$. This matrix corresponds to what is usually called the Seifert matrix for a sphere knot, so we define the  \emph{Seifert matrix} of a disk knot to be the integer matrix  $\theta=(\tau R)'$. Similarly,  the  linking matrix corresponding to $i_+$ is 
$L''(i_+p_*(\delta_j), \delta_i)=\sum_{k=1}^m R_{kj}\mu_{ik}$. Using the equality $L''(i_{-}p_*(\delta_i), \delta_j)=(-1)^{n+1}L''(i_{+}p_*(\delta_j), \delta_i)$ of \cite{GBF1}, we obtain that $\mu R=(-1)^{n+1}R'\tau'$. It is shown in \cite{GBF1} that $(-1)^{q+1}(R^{-1})'\tau R t-\tau'$ is a presentation matrix for $\cok(H_n(\td C;\Q)\to H_n(\td C, \td X;\Q))$ as module over $\Gamma=\Lambda\otimes_{\Z} \Q=\Q[\Z]$, while the matrix $\frac{1-t}{(R^{-1})'\tau t-(-1)^{q+1}\tau'R^{-1}}$ represents the Blanchfield pairing of this module. Both of these matrices are with respect to the natural  \emph{integral} bases within the rational modules. 
                     
Now, $L''(i_{-}p_*(\delta_i), \delta_j)-L(i_+p_*(\delta_i), \delta_j)=(\tau R)'_{ij}-(\mu R)'_{ij}=(\tau R)'_{ij}-((-1)^{n+1}(\tau R)')'_{ij}=\theta+(-1)^n\theta'$. But with respect to dual bases, this is also the matrix of $-R'$. So we have  $\theta+(-1)^n\theta'=-T$. 

Note that there is a correspondence between sphere knots and  disk knots whose boundary knots are trivial: Given such a disk knot, we can cone the boundary to obtain a \emph{locally-flat} sphere knot, and conversely, given a sphere knot, we can remove a ball neighborhood of any point on the knot to obtain a disk knot with trivial boundary. If we then consider a Seifert surface for such a disk knot whose boundary Seifert surface is the trivial disk Seifert surface for the boundary unknot, then the map $p_*$ will be an isomorphism (the identity with a proper choice of bases) and the Seifert matrix $\theta$ will be the same as the ordinary sphere knot Seifert matrix for the corresponding sphere knot. 

\vskip1cm
We close this introductory section by reminding the reader of some terminology.

 A  sphere knot $S^{n-2}\into S^{n}$  is called \emph{simple} if $\pi_i(X)\cong \pi_i(S^1)$ for $i\leq \frac{n-2}{2}$. By \cite{L66} this is as connected as the complement of a knot can be without the knot being trivial. We similarly define a disk knot  $D^{n-2}\into D^{n}$  to be \emph{simple} if $\pi_i(C)\cong \pi_i(S^1)$ for $i\leq \frac{n-2}{2}$. 

A pairing of modules $(\,,\,): A\otimes B\to C$ is called nondegenerate if $(a,b)=0$ for all $b\in B$ implies $a=0$ and if $(a,b)=0$ for all $a\in A$ implies $b=0$. We call the pairing nonsingular if $a\to (a, \cdot)$ is an isomorphism $A\to \Hom(B,C)$ and $b\to (\cdot , b)$ is an isomorphism $B\to \Hom(A,C)$. A rational matrix is nondegenerate and nonsingular if its determinant is not $0$. An integer matrix is considered nondegenerate if its determinant is nonzero and nonsingular if its determinant is $\pm 1$. 

\section{Disk knot cobordism}\label{S: cobordism}

Let $L: D^{n-2}\subset D^n$ be a disk knot. We define two types of cobordism between disk knots:

\begin{definition}
Two disk knots $L_0$, $L_1$ are cobordant if there exists a proper embedding $F: D^{n-2}\times I\into D^n\times I$ such that $F|  D^{n-2}\times {i}=L_i\times i$ for $i=0,1$ and $F| \bd D^{n-2} \times I$ is a cobordism of the boundary sphere knots $K_0$, $K_1$. 
\end{definition}

\begin{definition}
Two disk knots $L_0$, $L_1$ are cobordant rel boundary if the boundary knots of $L_0$ and $L_1$ are ambient isotopic and there exists a cobordism from $L_0$ to $L_1$ that restricts to this isotopy on $\bd D^{n-2}\times I$. 
\end{definition}

\textbf{N.B.} Due to the usual orientation switch of the total space from the bottom to the top of a cylinder, the embedding $L_1\times 1$ actually represents the knot $-L_1$, the mirror image of $L_1$. This will be the case, in particular, when we consider $L_1\times 1=-L_1$ as a submanifold of $S^{n}=\bd( D^n\times I)$. The orientation of the embedded knot is itself switched, of course, but this orientation usually plays no role in higher-dimensional knot theory so we omit further mention.

\begin{proposition} \label{P: even cob rel}
If $n$ is even, then any disk two knots $L_0, L_1: D^{n-2}\to D^n$ with isotopic boundary knots are  cobordant rel boundary. 
\end{proposition}

\begin{proof}
Given two such knots $L_0$, $L_1$, the maps $L_0\times i$ on $D_{n-2}\times i$ along with the isotopy $H$ connecting their boundary knots determines a sphere knot $S^{n-2}\to S^n=\bd(D^n\times I)$ by $(L_0\times 0)\cup H\cup (L_1\times 1):S^{n-2}=(D^{n-2}\times 0)\cup (S^{n-3})\times I\cup (D^{n-2}\times 1)\to S^{n}=(D^{n}\times 0)\cup (S^{n-1})\times I\cup (D^{n}\times 1)$. This is an even dimensional sphere knot and so it is null-bordant by Kervaire \cite{Ke}. Any such null-cobordism provides the desired cobordism of the disk knots. 

\end{proof}

\begin{proposition}\label{P: even cob}
If $n$ is even, then two disk knots $L_0, L_1: D^{n-2}\to D^n$ are cobordant if and only if their boundary knots are cobordant. 
\end{proposition}

\begin{proof}
The proof of the existence of a cobordism if the boundary knots are cobordant is the same as in the last proposition but connecting the boundary knots by their cobordism instead of the trace of an isotopy. The converse is immediate.
\end{proof}

This leaves the cases for which $n$ is odd. Note that in this case all boundary knots $K$ are cobordant since they will all be even dimensional \cite{Ke}.

\begin{proposition}\label{P: odd cob}
If $n$ is odd, then any disk two knots $L_0, L_1: D^{n-2}\to D^n$ are cobordant.
\end{proposition}
\begin{proof}
The boundary knots $K_0$, $K_1$ of $L_0$, $L_1$ will be even dimensional. As noted, all even dimensional knots are nullcobordant by Kervaire \cite{Ke}. Let us construct the cobordism $G$ of the boundary knots $K_0$ and $K_1$ as follows: Let $G|D^{n-2}\times [0,1/4]$ realize a null-cobordism of $K_0$. The union of $L_0$ with this nullcobordism givens a disk knot in $(D^{n}\times 0)\cup (\bd D^n\times [0,1/4])$ with unknotted boundary knot.  Let $J_0$ denote the sphere  knot obtained by filling in this unknotted boundary; we can think of obtaining $J_0$ by taking the cone pair on the boundary of the disk knot (which will be a locally-flat sphere knot since the boundary knot is trivial). Define $-J_1$ similarly by adjoining a null-cobordism on $\bd D^n\times [3/4,1]$ (recall that the embedding $L_1\times 1$ represents the disk knot $-L_1$, taking into account orientations on the cylinder as induced from the $0$ end). Now consider the knot $-(J_0 \# -J_1)=(-J_0)\#J_1$, where $\#$ represents knot sum. By removing neighborhoods of two points on the knot, we can think of this knot as a cobordism between two trivial knots $S^{n-3}\subset S^{n-1}$, and we can glue this cobordism into $\bd D^n\times [1/4,3/4]$, matching the ends since all unknots are ambient isotopic. So now we have constructed a cobordism from $K_0$ to $K_1$, and the knotted sphere in the boundary of $D^n\times I$ given by the union of $L_0$, $L_1\times 1=-L_1$, and the cobordism is $ J_0 \# \left((-J_0)\#J_1\right)\#(-J_1)=(J_0\#J_1)\#-(J_0\#J_1)$, which is null-cobordant. Again any null-cobordism now realizes the cobordism of disk knots. 
\end{proof}

So we are now reduced to the much more difficult consideration of cobordism of odd dimensional disk knots rel boundary. As seen in the preceding propositions and described in more detail below, the problem reduces to finding a null-cobordism of sphere knots composed of the union of $L_0$ and $L_1$. By \cite{L69} the cobordism class of a sphere knot $S^{2n-1}\into S^{2n+1}$, $n>1$, is determined by its Seifert matrix. So we are left with the problem of determining Seifert matrices for disk knots joined along their boundaries. Note that if $n=1$, the disk knot $D^1\into D^3$ has trivial boundary and so the problem of determining cobordisms rel boundary is in this case equivalent to the problem of classifying cobordisms of classical knots, which remains an unsolved problem. Hence we concentrate on the cases $m>1$ in which the disk knot and sphere knot theories are truly different (though closely related). 

We begin with a variety of algebraic preliminaries which generalize those of Levine in \cite{L69}. Levine defines a $2r\times 2r$ integer matrix to be null-cobordant if it is integrally  congruent to a matrix of the form $\begin{pmatrix}0 &N_1\\N_2&N_3\end{pmatrix}$, where each matrix $N_i$ is $r\times r$. Similarly, we will call a rational $2r\times 2r$ matrix $A$ rationally null-cobordant if it is rationally congruent to a matrix of this form. This is equivalent to saying that  $A$ is null-cobordant as a pairing of rational vector spaces $\Q^{2r}\times \Q^{2r}\to \Q$ ($x\times y\to x'Ay$), i.e. there exists an $r$-dimensional subspace of $\Q^{2r}$ on which the restriction of the pairing is $0$. 

\begin{lemma}\label{L: int is rat}
Let $A$ be a $2r\times 2r$ integral matrix. Then $A$ is null-cobordant if and only if it is rationally null-cobordant.
\end{lemma}

\begin{proof}
If $A$ is a $2r\times 2r$ integral null-cobordant matrix, then there is a rank $r$ direct summand $F$ of $\Z^{2r}$ on which $A$ restricts to the $0$ bilinear form. Hence this matrix is also rationally null-cobordant, restricting to the $0$ form on $F\otimes \Q$. 

Conversely, suppose that $A$ is rationally null-cobordant so that there is an $r$-dimensional $\Q$ subspace $V$ of $\Q^{2r}\cong \Z^{2r}\otimes \Q$ on which $A$ restricts to the $0$ bilinear form. Let $L$ be the lattice $\Z^{2r}\cap V$. This is a free abelian subgroup of $\Z^{2r}$, in fact a direct summand since any element of $\Z^{2r}$ that has a scalar multiple in $L$ must also be in $L$. $L$ must have rank at least $r$, since given $r$ linear independent rational vectors in $V$, there are integral multiples of these vectors that lie in $L$ (by clearing denominators of the coordinates), and these scalar multiples remain linearly independent over $\Q$ and hence over $\Z$. So $A$ is the $0$ form on a free abelian group of rank $\geq r$ that is a direct summand of $\Z^{2r}$.
\end{proof}

\begin{corollary}
If $A$ is a rationally null-cobordant matrix obtained from an integral Seifert matrix of a $2n+1$ sphere knot $K$, $n>1$, by a rationally unimodular congruence, then the knot is null-cobordant.
\end{corollary}

\begin{proof}
Clearly any matrix rationally unimodularly congruent to a rationally null-cobordant matrix is also rationally null-cobordant. So by the preceding lemma, the integral Seifert matrix of $K$ is null-cobordant, and the result follows from the Main Theorem of \cite{L69}. 
\end{proof}

\begin{corollary}\label{C: rational nullcobordance}
Let $A$ be the matrix of the rational Seifert pairing of a $2n+1$ sphere knot, $n>1$, with respect to some Seifert surface $F$ and some rational basis of $H_n(F;\Q)$. Then $A$ is rationally null-cobordant if and only if the knot is null-cobordant. 
\end{corollary}
\begin{proof}
There is a rationally unimodular change of basis that will take the rational basis for $H_n(F;\Q)$, with respect to which $A$ is defined, to an integer basis of the group $H_n(F;\Z)$. In this basis, we obtain an integral Seifert matrix $B$ for the knot. If the knot is null-cobordant, there will be a rank $r$ summand of $H_n(F;\Z)$ on which the form determined by $B$ is $0$. The tensor product of this summand with $\Q$ gives a self-annihilating $r$-dimensional $\Q$ vector subspace of $H_n(F;\Q)$. Under any rationally unimodular change of basis, this subspace (or more precisely, it's image under the change of basis) will remain self-annihilating. In particular, $A$ will be rationally null-cobordant.

Conversely, if $A$ is rationally null-cobordant, then the knot is null-cobordant by the previous corollary.
\end{proof}

\begin{definition}
Two square rational matrices $A$ and $B$ are rationally cobordant if $A\boxplus -B$ is rationally null-cobordant, where $\boxplus$ denotes the block sum of matrices $A\boxplus -B=\begin{pmatrix}A &0\\0&-B\end{pmatrix}$.
\end{definition}

\begin{corollary}\label{C: rat gives int}
Two \emph{integral matrices} are integrally cobordant if and only if they are rationally cobordant. \end{corollary}
\begin{proof}
This is an immediate consequence of Lemma \ref{L: int is rat}.
\end{proof}

\textbf{N.B.} Even though we will be concerned with rational cobordism class, the term \emph{Seifert matrix} will always refer to the \emph{integral} Seifert matrix defined in Section \ref{S: basics} \emph{unless explicitly stated otherwise}.

\begin{lemma}
Suppose $A$ and $N$ are square matrices of rational numbers, that $N$ and $A\boxplus N$ are rationally null-cobordant, and that some rational linear combination $\lambda N+\mu N'$ has non-zero determinant. Then $A$ is rationally null-cobordant. 
\end{lemma}
\begin{proof}
The proof is the same as that of Levine's for integral null-cobordism \cite[Lemma 1]{L69} replacing $\Z$ with $\Q$ in all steps. 
\end{proof}

For Seifert matrices of sphere knots $S^{2n-1}\subset S^{2n+1}$ or disk knots $D^{2n-1}\subset D^{2n+1}$, these conditions will be satisfied with $\lambda=1$, $\mu=(-1)^n$. For sphere knots, this is well-known (see \cite{L66} or \cite{L69}). For disk knots, this can be concluded similarly from the fact that the Alexander polynomials of disk knots are non-zero when evaluated at $1$ (see \cite{GBF} or \cite{GBF1} for details). 

\begin{corollary}
For fixed rational $\lambda$ and $\mu$, the set of rational cobordism classes of square rational matrices $A$ satisfying $\det(\lambda A +\mu A')\neq 0$ is an abelian group under block sum, the inverse of the class represented by a matrix $A$ being the class represented by $-A$. 
\end{corollary} 
\begin{proof}
Again, this corollary follows from the  lemma as in \cite[\S 3]{L69} by replacing integral statements with rational ones. 
\end{proof}

\begin{proposition}\label{L: s-equivalence cobordism}
Let $A$, $\hat A$ be Seifert matrices for the disk knot $L$. Then  $A$ and $\hat A$ are integrally cobordant.   
\end{proposition}
\begin{proof}
Let $V$ and $\hat V$ be Seifert surfaces with respect to which $A$ and  $\hat A$ are the integral Seifert matrices. Then it follows from the results of Section \ref{S: seifert surfaces}, below, that $A$ and $B$ are related by a sequence of rational congruences and enlargements or reductions of the form
\begin{equation*}
M
\leftrightarrow M'=
\begin{pmatrix}
M&0&\eta\\
0&0&x\\
\xi &x'&y
\end{pmatrix},
\end{equation*}
where $M$ is a matrix, $\eta$ is a column vector, $\xi$ is a row vector, $x$, $x'$, and $y$ are integers, and all ``$0$''s represent the necessary $0$ entries to make this matrix square. Also, one of $x, x'$ is $0$ while the other is non-zero. So it suffices to show that $-M\boxplus M'$ is rationally nullcobordant. If $M$ is a $k\times k$ matrix, let $I_j$ be the $j\times j$ identity matrix, and let $P=\begin{pmatrix}I_k& I_k &0\\0&I_k&0\\0&0&I_2\end{pmatrix}$. Then $P'(-M\boxplus M')P=\begin{pmatrix}
-M&-M&0&0\\
-M&0&0&\eta\\
0&0&0&x\\
0&\xi &x'&y
\end{pmatrix}$ contains a $k+1\times k+1$ dimensional $0$ matrix block symmetric with respect to the diagonal, so it is rationally null-cobordant. The integral cobordism is then implied by Corollary \ref{C: rat gives int}.
\end{proof}

Let $L_0$ and $L_1$ be two $2n+1$ disk knots with the same boundary knot $K$. Then these knots will be cobordant rel boundary if and only if the knot $\mf{K}=L_0\cup_{K} -L_1$ is null-cobordant as a sphere knot. So we must examine its Seifert matrix. Let $V_0$ and $V_1$ be Seifert surfaces for $L_0$ and $L_1$  with boundary Seifert surfaces $F_0$ and $F_1$ for $K$ (see \cite{GBF1}). Then there is a cobordism $\Upsilon$ of Seifert surfaces from $F_0$ to $F_1$ with boundary the union of $F_0$, $-F_1$, and the trace of an isotopy of $K$ \cite[\S 3]{L70}. Then we can form a Seifert surface for $\mf K$ by $W=V_0\cup_{F_0} \Upsilon\cup_{-F_1} -V_1$. Since the union of $L_0$ with the trace of an isotopy of its boundary is isotopic to $L_0$, we will simplify notation by combining $V_0$ and $\Upsilon$ to form a new $V_0$. So we can consider $W$ to be composed of  Seifert surface $V_0$ and $-V_1$ for $L_0$ and $-L_1$, joined along a single Seifert surface $F$ for $K$.  

In what follows, we use the isomorphism of the groups $H_*(V_1)\cong H_*(-V_1)$ to simplify the notation.

We consider the Mayer-Vietoris sequence
\begin{equation*}
\begin{CD}
@>>> & H_n(F) &@>j>>& H_n(V_0)\oplus H_n(V_1) &@>\rho>>& H_n(W) &@>\bd >>& H_{n-1}(F) &@>j'>>.
\end{CD}
\end{equation*}
We are first interested in computing ranks of free abelian subgroups, so we can consider homology groups with \emph{rational} coefficients (though we omit them from the notation for clarity). Then there is a splitting $H_n(W)\cong \im(\bd)\oplus \cok(j)$. 

Now from the \emph{rational} long exact sequences of the pairs $(V_s, F)$, $s=0,1$:
\begin{equation}\label{E: Seifert pair}
\begin{CD}
@>>> & H_n(F) &@>i_s>>& H_n(V_s) &@>p_s>>& H_n(V_s,F) &@>\bd_s >>& H_{n-1}(F) &@>i_s'>>,
\end{CD}
\end{equation}
$H_n(V_s)\cong \cok(i_s)\oplus \im(i_s)$, and, furthermore, $\im(j)\cap H_n(V_s)\subset \im(i_s)$, so we can write $\cok(j)\cong \frac{\cok(i_0)\oplus \cok(i_1)\oplus \im(i_0)\oplus \im(i_1)}{\im(j)}\cong \cok(i_0)\oplus \cok(i_1)\oplus \frac{\im(i_0)\oplus \im(i_1)}{\im(j)}$. 

Now $\cok(i_s)$ is the group on which the Seifert matrix of $L_s$ is defined. We need to study the other summands $\frac{\im(i_0)\oplus \im(i_1)}{\im(j)}$ and $\im(\bd)$ of $H_n(W)$. We claim that these two summands have the same dimension. 

Let $|G|$ stand for the dimension of the vector space $G$. Suppose that $|H_n(F)|=m$ and $H_n(V_s)=M_s$. Then $|\im(i_s)|=|\coim(i_s)|=m-|\ker(i_s)|$, and $|\im(j)|= m-|\ker(j)|=m-|\ker(i_0)\cap \ker(i_1)|$. So, 
\begin{align*}
\left|\frac{\im(i_0)\oplus \im(i_1)}{\im(j)}\right|&= m-|\ker(i_0)|+ m-|\ker(i_1)|-( m-|\ker(i_0)\cap \ker(i_1)|)\\
&=m+|\ker(i_0)\cap \ker(i_1)| -|\ker(i_0)|-|\ker(i_1)|.
\end{align*}

Now, since $F$ is a $2n-1$, $n>1$, manifold with sphere boundary, and since $V_s$ is a $2n$-manifold whose boundary is the union of $F$ with a disk, Poincar\'e duality holds, and, in particular,  $|H_n(F)|=|H_{n-1}(F)|=m$ and $|H_n(V_s)|=|H_n(V_s,F)|+M_s$. Let us fix a basis of $H_{n-1}(F)$ and use the standard orthonormal inner product with respect to this basis to identify $H_{n-1}(F)$ with $\Hom(H_{n-1}(F);\Q)\cong H^{n-1}(F;\Q)\cong H_n(F;\Q)$. Consider now $\ker(i_s)$. Under this identification, via Poincar\'e duality, $\ker(i_s)=(\ker(i_s'))^{\bot}$. Indeed, if $x\in \ker(i_s)$ and $y\in \ker(i_s')=\im(\bd_s)$ so $y=\bd_s z$, then on the intersection pairing, we have $S_F(x,y)=S_V(i_s(x),z)=0$; so the identification takes $\ker(i_s)$ into $(\ker(i_s'))^{\bot}$. But also we have $|\ker(i_s)|=m-|\im(i_s)|=m-|\ker(p_s)|=m-(M_s-|\im(p_s)|)=m-(M_s-|\ker(\bd_s)|)=m-(M_s-(M_s-|\im(\bd_s)|))=m-|\ker(i'_s)|$, so $\ker(i_s)=(\ker(i_s'))^{\bot}$. Then we compute

\begin{align*}
\left|\frac{\im(i_0)\oplus \im(i_0)}{\im(j)}\right|
&=m+|\ker(i_0)\cap \ker(i_1)| -|\ker(i_0)|-|\ker(i_1)|\\
&=m+|(\ker(i'_0))^{\bot}\cap (\ker(i'_1))^{\bot}|- |(\ker(i'_0))^{\bot}|-|(\ker(i'_1))^{\bot}|\\
&=m-|(\ker(i'_0))^{\bot}+(\ker(i'_1))^{\bot}|\\
&=|\ker(i'_0)\cap\ker(i'_1)|\\
&=|\ker j'|\\
&=|\im(\bd)|.
\end{align*}

Here the fourth equality uses that, in a vector space $X$ with subspaces $Y$ and $Z$, $(Y^{\bot}+Z^{\bot})^{\bot}=Y\cap Z$: if $x\in Y\cap Z$, then $\langle x, a\rangle=0$ if $a$ in $Y^{\bot}$ or $Z^{\bot}$, so this is $0$ if $a\in (Y^{\bot}+Z^{\bot})$, so $Y\cap Z\subset (Y^{\bot}+Z^{\bot})^{\bot}$. Conversely, if $x\in (Y^{\bot}+Z^{\bot})^{\bot}$, then $\langle x,a\rangle=0$ for any a in $Y^{\bot}+Z^{\bot}$ and, in particular, any $a$ in either $Y^{\bot}$ or $Z^{\bot}$. So $a\in Y^{\bot\bot}=Y$ and similarly $a\in Z$. So $(Y^{\bot}+Z^{\bot})^{\bot}\subset Y\cap Z$. Hence, $((\ker(i'_0))^{\bot}+(\ker(i'_1))^{\bot})^{\bot}=\ker(i'_0)\cap\ker(i'_1)$, so these spaces have complementary dimensions in $H_{n-1}(F)$.

So once again, with rational coefficients, we can write $H_n(W)\cong \cok(i_0)\oplus \cok(i_1)\oplus \frac{\im(i_0)\oplus \im(i_1)}{\im(j)} \oplus \im(\bd)$, where the last two summands have the same dimension. Let us denote $U=\frac{\im(i_0)\oplus \im(i_1)}{\im(j)}$

We next observe that the Seifert form is $0$ when restricted to $U\times (\cok(i_0)\oplus \cok(i_1))$ or $(\cok(i_0)\oplus \cok(i_1))\times U$. This is true because any element of $\cok(i_s)$ can  represented by a cycle lying entirely in the interior of   $V_s$ and hence of  $D^{2n+1}\times s$  in the cobordism, and the same is true for any translate along a normal vector to the Seifert surface. Also, we can then find a chain in $D^{2n+1}\times s$ whose boundary is the push in the bicollar  of our cycle.  Meanwhile, any element of $U$ can be represented by a cycle that lies in  $\bd D^{2n+1}\times I$, and the same for its translates along the bicollar, and a choice of chain it bounds in $\bd D^{2n+1}\times I$. So then clearly the linking numbers of any such cycles must be $0$.

At last we can prove the following proposition.

\begin{theorem}\label{T: combo knot}
Let $A_0$ and $A_1$ be Seifert matrices for $2n+1$ disk knots $L_0$ and $L_1$, $n>1$, with the same boundary knot. Then the Seifert matrix of the sphere knot $L_0\cup_{\bd}-L_1$ is integrally cobordant to  $A_0\boxplus  -A_1$.
\end{theorem}

\begin{proof}
Let $V_0$, $V_1$, $F$, and $W$ be as above. Let $B_0$ and $B_1$ be the Seifert matrices of $L_0$ and $L_1$ corresponding to these Seifert surfaces. Then by Proposition \ref{L: s-equivalence cobordism}, $B_0$ and $B_1$ are rationally cobordant to $A_0$ and $A_1$, respectively.

Now we consider the Seifert matrix $M$ determined by $W$ and show that is is rationally cobordant to $A_0\boxplus  -A_1$, which will suffice to prove the theorem. 

We know that $H_n(W;\Q)\cong  \cok(i_0)\oplus \cok(i_1)\oplus \frac{\im(i_0)\oplus \im(i_1)}{\im(j)} \oplus \im(\bd)$, and the Seifert pairings on $\cok(i_0)$ and $\cok(i_1)$ must restrict to $B_0$ and $-B_1$ by definition (the negative is due to the reverse of orientation by considering $L_1$ in $D^{2n+1}\times 1\subset D^{2n+1}\times I$). Furthermore, these subspaces are orthogonal under the Seifert pairing since elements of $\cok(i_0)$ are represented by chains in $V_0\subset D^{2n+1}\times 0 \subset \bd (D^{2n+1}\times I)$, while elements of  $\cok(i_1)$ are  represented by chains in $V_1\subset D^{2n+1}\times 1 \subset \bd( D^{2n+1}\times I)$, so these chains cannot link in $\bd (D^{2n+1}\times I)=S^{2n+1}$. Similarly, elements in $\frac{\im(i_0)\oplus \im(i_0)}{\im(j)}$ can be represented by chains in $F$ that can be pushed into either $V_0$ or $V_1$ and so these do not link with each other or elements of $\cok(i_0)$ and $\cok(i_1)$. Thus $M$ must have the form (up to rational change of basis and hence rational cobordism)
\begin{equation*}
M=
\begin{pmatrix}
B_0 & 0       & 0 & X_1 \\
0     & -B_1  & 0 & X_2 \\
0     & 0       &0   & X_3\\
X_4 & X_5  & X_6 & X_7
\end{pmatrix}
\end{equation*}
for some matrices $X_i$. Note that the diagonal blocks are all square and that the last two diagonal blocks have the same size by the above dimension calculations.  This is a generalization of the kind of elementary enlargement that we considered in Proposition \ref{L: s-equivalence cobordism}. Set $\mf P=\begin{pmatrix}I_r& I_r &0\\0&I_r&0\\0&0&I_{2s}\end{pmatrix}$, where $r=|\cok(i_0)\oplus \cok(i_1)|$ and $s=|\frac{\im(i_0)\oplus \im(i_0)}{\im(j)} |=|\im(\bd)|$. Then 
$$\mf P'(-B_0\boxplus B_1 \boxplus M)\mf P=
\begin{pmatrix}
-B_0&  0  &-B_0      & 0       &0  & 0 \\
0&  B_1  &0      & B_1       &0  & 0 \\
-B_0&  0  &0     & 0       & 0 & X_1 \\
0&  B_1  &0     & 0      & 0 & X_2 \\
0&  0  &0     & 0       &0   & X_3\\
0&0&X_4 & X_5  & X_6 & X_7
\end{pmatrix}
$$ contains an $r+s \times r+s$ trivial submatrix symmetric about the diagonal. Hence it is rationally null-cobordant and $M$ is rationally cobordant to $B_0\boxplus -B_1$, which in turn is rationally cobordant to $A_0\boxplus -A_1$ using Proposition \ref{L: s-equivalence cobordism}. The rational cobordisms become integral cobordisms by Corollary \ref{C: rat gives int}.
\end{proof}

The following theorem now follows immediately. 
\begin{theorem}\label{T: main cobordism}
Let $A_0$ and $A_1$ be Seifert matrices for disk knots $L_0, L_1: D^{2n-1}\into D^{2n+1}$ with the same boundary knot. $L_0$ and $L_1$ are cobordant rel boundary if and only if $A_0$ and $ A_1$ are cobordant.
\end{theorem}
\begin{proof}
If the matrices are cobordant, then the integral Seifert matrix for $L_0\cup_{\bd}-L_1$, which is rationally cobordant to  $A_0\boxplus  -A_1$, is rationally nullcobordant, hence integrally null-cobordant. Thus $L_0\cup_{\bd}-L_1$ is slice and the slicing disk provides the desired cobordism. Conversely, if $A_0$ and $ A_1$ are not cobordant, then $A_0\boxplus  -A_1$ is not integrally nullcobordant, so there can be no such slicing disk to provide the cobordism. 
\end{proof}

\begin{corollary}\label{C: Alex}
Suppose that $L_0$ and $L_1$ are disk knots $D^{2n-1}\subset D^{2n+1}$, $n>1$, such that $\bd L_0=\bd L_1=K$. Then a necessary condition for $L_0$ and $L_1$ to be cobordant rel boundary is that the product of the middle-dimensional Alexander polynomials $c^{L_0}_{n}(t)c_{n}^{L_1}(t)$ be similar in $\Q[t,t^{-1}]$ to a polynomial of the form $p(t)p(t^{-1})$.
\end{corollary}
\begin{proof}
By \cite[\S 3.6]{GBF1} or \cite[\S 3.6]{GBF} and the calculations in Section \ref{S: basics} above, $c^{L_i}_n(t)$ is in the similarity class in $\Q[t,t^{-1}]$ of the determinant of $(A_i+(-1)^{n}A'_i)^{-1}(A_it+(-1)^{n}A_i')$, where $A_i$ is the Seifert matrix of $L_i$, $i=0,1$. We know that if $L_0$ and $L_1$ are cobordant rel boundary, then $B=A_0\boxplus -A_1$ is rationally nullcobordant. It follows then as in \cite[\S 15]{L69} that the determinant of $Bt+(-1)^nB'$ is similar to  $\bar p(t)\bar p(t^{-1})$ for some polynomial $\bar p$. But clearly the determinant of $Bt+(-1)^nB'$ is equal to $\pm$ the product of the determinants of $(A_it+(-1)^{n}A_i')$, $i=0,1$. So $c^{L_0}_{n}(t)c_{n}^{L_1}(t)\sim \frac{\bar p(t)\bar p(t^{-1})}{\det(A_0+(-1)^{n}A'_0)\det(A_1+(-1)^{n}A'_1)}$. The claim now follows since $\frac{1}{\det(A_0+(-1)^{n}A'_0)\det(A_1+(-1)^{n}A'_1)}$ is a unit in $\Q[t,t^{-1}]$.
\end{proof}

Since $L_0\cup_{K} -L_1$ is a sphere knot $S^{2n-1}\subset S^{2n+1}$, there is a basis for which its integral  Seifert matrix $A$ is a matrix of integers such that $A+(-1)^nA'$ is integrally unimodular. Thus each possible obstruction matrix $A_0\boxplus -A_1$ must be rationally cobordant to such a matrix. We can also state the following converse:

\begin{theorem}
Let $A$ be a matrix of integers such that $A+(-1)^nA'$ is integrally unimodular, and let $L_0$ be a disk knot $D^{2n-1}\subset D^{2n+1}$ with Seifert matrix $A_0$. Then there is a disk knot $L_1$ with the same boundary knot  as $L_0$ and such that the obstruction Seifert matrix $A_0\boxplus -A_1$ to $L_0$ and $L_1$ being cobordant rel boundary is cobordant to $A$. 
\end{theorem}
\begin{proof}
By \cite{L69}, there is a sphere knot $\mc{K}: S^{2n-1}\subset S^{2n+1}$ with Seifert matrix $-A$. Let $L_1$ be the knot $L_0\#\mc{K}$, the knot sum taken away from the boundary. Then $L_1$ has Seifert matrix $A_0\boxplus -A$, and $L_0\cup_K -L_1$ has Seifert matrix rationally cobordant to $A_0 \boxplus -A_0 \boxplus A$, which is rationally cobordant to $A$, hence integrally cobordant to $A$ by Corollary \ref{C: rat gives int}. 
\end{proof}

Similarly, we can show the following:

\begin{theorem}\label{T: there is a sphere}
Let $L_0$ and $L_1$ be two disk knots $D^{2n-1}\subset D^{2n+1}$, $n>1$, with common boundary $K$. Then there exists a sphere knot $\mc{K}: S^{2n-1}\subset S^{2n+1}$ such that $L_0$ is cobordant to the knot sum (away from the boundary) $L_1\# \mc {K}$. 
\end{theorem}
\begin{proof}
Let $A_0$ and $A_1$ be the Seifert matrices for $L_0$ and $L_1$. Then as above $A_0\boxplus -A_1$ is rationally cobordant to an integral matrix $B$ such that $B+(-1)^nB'$ is integrally unimodular. Let $\mc{K}$ be a sphere knot with Seifert matrix $B$, which exists by \cite{L69}. Then $L_1\# \mc{K}$ has Seifert matrix $A_1\boxplus B$, and  $L_0\cup_K -(L_1\#\mc{K})$ has Seifert matrix $A_0\boxplus -A_1\boxplus -B$, which is rationally null-cobordant. So $A_0$ is rationally cobordant to $A_1\boxplus B$, and the theorem now follows from Theorem \ref{T: main cobordism} and Corollary \ref{C: rat gives int}.
\end{proof}

\begin{theorem}
Let $K: S^{2n-2}\into S^{2n}$, $n>1$ be a sphere knot. Then there is a disk knot $L: D^{2n-1}\into D^{2n_1}$ such that $\bd L=K$ and $\pi_i(D^{2n+1}-D^{2n-1})\cong \pi_i(S^1)$ for $i<n$. 
\end{theorem}
\begin{proof}
By Kervaire \cite[Thm. III.6]{Ke}, there exists some disk knot whose boundary is $K$ (all even dimensional knots are null-cobordant). We show that in fact Kervaire's construction gives us a knot of the desired type. The argument in Kervaire's theorem proceeds as follows (modifying the notation slightly to coincide with our own): Let $F$ be a Seifert surface for $K$. Then it is possible to construct a manifold $V^{2n}$ and to embed it into $D^{2n+1}$ such that $V\cap S^{2n-1}=F$ and $\bd V=F\cup D^{2n+1}$. This manifold $V$ will be a Seifert surface for $L$, and it is obtained from $F$ by adding handles of core dimension $\leq n$ to $F\times I$, in order of increasing dimension, to successively kill the homotopy groups of $F$ by surgery. In particular then, after the addition of the $2$-handles to $F\times I$, we obtain a simple connected manifold as the trace of the surgery, and ultimately $H_{2n-i}(V,F)=0$ for $i<n$ because there are no handles of core dimension $>n$ added. Then $H_{2n-i}(V,F)\cong H^{i}(V)$ for $i\geq 1$, so $H^i(V)=0$ for $1\leq i<n$, which implies that $H_i(V)=0$ for for $1\leq i<n$. 

It now follows that $D^{2n+1}-V$ is simply-connected by the van Kampen theorem: by pushing along the bicollar of $V$, we can thicken $V$ to a homotopy equivalent $2n+1$ manifold  whose common boundary with the closure of its complement in $D^{2n+1}$ is the union of two copies of $V$ glued along $L$ (see \cite{L65}). It then follows from the van Kampen theorem that $D^{2n+1}-V$ must be simply-connected, and from Alexander duality for a ball that $H_i(D^{2n+1}-V)=0$ for $0<i<n$ (see \cite[Prop. 3.3]{GBF1} and note that these arguments extend to integer coefficients). 

Now, using the usual cut-and-past construction of the infinite cyclic cover of $D^{2n+1}-L$ (see \cite{L66}), another inductive application of the van Kampen theorem shows now that the infinite cyclic cover of $D^{2n+1}-L$ is simply connected, and the Mayer-Vietoris theorem shows that its homology is trivial in dimensions $<n$. So this cover is $n-1$-connected, and it follows that the homotopy groups $\pi_i(D^{2n+1}-L)$ vanish for $1<i<n$ and that $\pi_1(D^{2n+1}-L))\cong \Z$. 
\end{proof}

\begin{theorem}\label{T: cob to simple}
Given any disk knot $L:D^{2n-1}\subset D^{2n+1}$, $n>1$, $L$ is cobordant rel boundary to a disk knot $L_1$ such that $\pi_i(D^{2n+1}-L_1)\cong \pi_i(S^1)$ for $i<n$.
\end{theorem}
\begin{proof}
First assume $n>2$.
By the preceding theorem, there exists a disk knot $L_0$ whose boundary agrees with that of $L$ and which satisfies the require homotopy conditions. Let $A$ and $A_0$ be the respective Seifert matrices of $L$ and $L_0$. Then we know that the matrix $A\boxplus -A_0$ is rationally cobordant to an integral matrix $B$  such that the determinant of $B+(-1)^nB'$ is integrally unimodular since this is true for the integral Seifert matrix of the sphere knot $L\cup -L_0$. By Levine \cite{L69}, there exists a sphere knot $\mc{K}: S^{2n-1}\subset S^{2n+1}$ whose Seifert matrix is $B$ and such that  $\pi_i(S^{2n+1}-\mc{K})\cong \pi_i(S^1)$ for $i<n$. Let $L_1$ be the knot sum $L_0\#\mc{K}$ along the interior. Then $L_1$ satisfies the desired homotopy properties and has Seifert matrix $A_0\boxplus B$, which we know is cobordant to $A$ since $A\boxplus -A_0\boxplus -B$ is rationally cobordant to $B\boxplus -B$, which is null-cobordant. By Theorem \ref{T: main cobordism}, $L$ and $L_1$ are cobordant rel boundary.

If $n=2$, then \cite{L69} provides a $\mc{K}$ only if $B+B'$ has signature a multiple of $16$. But since $L\cup-L_0$ is a knot $S^2\subset S^4$, its Seifert matrices will all satisfy this property (again see \cite{L69}), hence so will $A\boxplus -A_0$ since signature is a matrix cobordism invariant (integrally and rationally). Thus the argument of the preceding paragraph applies again.
 \end{proof}

\section{Realization of cobordism classes}\label{S: realize}

Up to this point we have shown that two odd-dimensional disk knots are cobordant rel boundary if and only if their Seifert matrices are  cobordant. This leads to the natural question: what cobordism classes of matrices can be realized as the Seifert matrices of disk knots? We first demonstrate that we are truly dealing with a wider variety of objects than just Seifert matrices of sphere knots:

\begin{proposition}
There exist Seifert matrices for disk knots that are not  cobordant to Seifert matrices of sphere knots. In particular this implies that there are Seifert matrices for disk knots that are not  cobordant to any integer matrix $A$ such that $A+(-1)^nA'$ is integrally unimodular. 
\end{proposition}
\begin{proof}
Suppose, to the contrary, that every disk knot Seifert matrix is cobordant to some sphere knot Seifert matrix. Let us then fix a disk knot $L:D^{2n-1}\subset D^{2n+1}$, $n$ even, $n>2$, with some Seifert surface and with Seifert matrix $B$. By assumption, $B$ is  cobordant to a Seifert matrix $C$ of some sphere knot; this implies that $B$ must have an even number of rows and columns, since this must be true of $C$ (see, e.g., \cite[p. 178]{T73}). By \cite{L69}, there exists a sphere knot $K$ with a Seifert surface that realizes the Seifert matrix $-C$.  Therefore, the  knot sum $K\#L$ with Seifert surface given as the boundary connected sum of the Seifert surfaces of $K$ and $L$ will yield the  null-cobordant Seifert matrix $A=B\boxplus -C$. It then follows as in \cite[\S15]{L69}  that the determinant of $tA+A'$ is the product of $\pm$ a power of $t$ with a Laurent polynomial of the form $p(t)p(t^{-1})$. In particular, $|-A+A'|$ is $\pm$ a  square.

Now,  by \cite[\S 3]{GBF1} and the calculations of Section \ref{S: seifert surfaces}, below, the middle dimensional Alexander polynomial $c_n(t)$ of a disk knot, $n$ even,  is given, up to similarity, by the determinant of $(A+A')^{-1}(At+A')$, which, with our current assumptions, must thus be of the form $\frac{p(t)p(t^{-1})}{(p(1))^2}$ (up to similarity). In particular, we see that the value $c_n(-1)$ associated to $K\#L$ must be a square. But we also know that the Alexander polynomial of a direct sum is the product of the polynomials so that $c_n^{K\#L}\sim c_n^Kc_n^L$, where $\sim$ denotes similarity and we have labeled the polynomials with their knots in the obvious way. But $c_n^K(-1)$ must be $\pm$ a square since $K$ is a sphere knot \cite{L66}. So it would follow that $c_n^L(-1)$ must also always be a square. However, this contradicts the calculations in \cite[\S3.64]{GBF1} which demonstrate that any odd number can be realized as $c_n^L(-1)$ for some $L$ of our fixed dimension. 

Hence we have demonstrated, at least for $n$ even, that there must exist disk knot Seifert matrices that are not  cobordant to sphere knot Seifert matrices. 
\end{proof}

However, we do have the following proposition:

\begin{proposition}
Suppose that $A$ is the Seifert matrix of a disk knot $L: D^{2n-1}\into D^{2n+1}$ with boundary knot $K$. Then $B$ is in the cobordism class of  a Seifert matrix of a disk knot $L'$ with the same boundary $K$ if and only if $A\boxplus -B$ is cobordant to the Seifert matrix of a sphere knot $\mf{K}: S^{2n-1}\into S^{2n+1}$.
\end{proposition}
\begin{proof}
If $L$ and $L'$ are disk knots with the same boundary sphere knot $K$ and respective Seifert matrices in the cobordism classes of $A$ and $B$, then we can form the knot $\mf{K}=L\cup_K -L'$ by gluing $L$ and $-L$ together, identifying the boundaries $K$ and $-K$. By Theorem \ref{T: combo knot}, $A\boxplus -B$  is cobordant to the Seifert matrix of the sphere knot $\mf{K}$. 

Conversely, suppose that $A\boxplus -B$ is cobordant to the Seifert matrix of some sphere knot $\mf{K}$. Then $-A\boxplus B$ will be the Seifert matrix of $-\mf{K}$. Form $L'=L\#-\mf{K}$, the internal knot sum. The Seifert matrix of this $L'$ will be  the sum of $A$ with the Seifert matrix of $-\mf{K}$ and hence will be cobordant to $A\boxplus -A\boxplus B$, which is cobordant to $B$.  
\end{proof}

This proposition tells us how to recognize rational cobordism classes of Seifert matrices for disk knots with a given sphere knot \emph{provided that we already have a cobordism class of Seifert matrices with which to compare}. This is a nice start, but we would like to find a way to determine which cobordism classes are realizable starting only with information about the boundary knot. It will turn out that the crucial datum is supplied by the isometry class of the Farber-Levine torsion pairing on $T_{n-1}(\td X)$, the $\Z$-torsion subgroup of $H_{n-1}(\td X)$, so long as this group has no $2$-torsion.

Let us begin by examining further the necessary conditions for a matrix  $\theta$ to be a Seifert matrix for a disk knot. 
We  know from Section \ref{S: basics} that if we choose dual bases of $H_n(V)$ and $H_n(V,\bd V)$, then the matrix of the map $p_*: \bar E \to \td E$ will also represent the transpose of the self intersection pairing on $H_n(V)$. By the computations in that section, we have $R=-\theta'+(-1)^{n+1}\theta$. To emphasize this dependence, we will sometimes write $R=R_{\theta}$. Note also that since, with these bases, $\theta=(\tau R)' $ and $R$ is invertible, $\tau$ is also determined by $\theta$ as $\tau_{\theta}=\theta'R_{\theta}^{-1}=\theta'(-\theta'+(-1)^{n+1}\theta)^{-1}$ (unfortunately, we can't simplify this further since in general $\theta$ won't be invertible). 

Another necessary conditions is that  $\tau_{\theta}=\theta'R_{\theta}^{-1}=\theta'(-\theta'+(-1)^{n+1}\theta)^{-1}$ must be integral since $\tau'$ is the matrix of $L'' (i_{-*}(\alpha_j),\delta_n)$. Also, we must have $(R_{\theta}^{-1})'\tau_{\theta} R_{\theta}=(-\theta+(-1)^{n+1}\theta')^{-1}\theta'$ integral, since this is, up to sign, the matrix $\mu_{\theta}$, where $\mu$ is the matrix of $L'' (i_{+*}(\alpha_j),\delta_n)$. 

We note one implication of these requirements:

\begin{proposition}
Let $\theta$ be the Seifert matrix of a disk knot $D^{2n-1}\subset D^{2n+1}$. If $n$ is odd, or if $n$ is even and $\det(R_{\theta})\neq 0$ mod $2$, then  $\theta$ must be even dimensional (have an even number of rows and columns).
\end{proposition}
\begin{proof}
If $n$ is odd, then $R_{\theta}=-\theta'+(-1)^{n+1}\theta$ is skew-symmetric. But $R_{\theta}$ is nondegenerate, so it must have even dimension.

Next, suppose that $n$ is even. By \cite[\S3.6]{GBF1}, the \emph{integral} Alexander polynomial $c_n(t)$ of $L$ is the determinant of $(R^{-1})'\tau Rt+\tau'$  up to similarity in $\Lambda$. This equals $\frac{1}{\det(R)}\det(\tau R t+R'\tau')$. Now,  again by  \cite[\S3.6]{GBF1}, $c_n(t)\sim c_n(t^{-1})$, but this really follows just from the symmetry of the presentation and so holds for the determinant of any matrix of the form  $(R^{-1})'\tau Rt+\tau'$. The same is also true of $\det(\tau R t+R'\tau')$, which we will call $d_n(t)$ (so $c_n(t)\sim d_n(t)/\det(R)$). By multiplying by a power of $t$, we can assume that $d_n(t)=\sum_{i=0}^m a_it^i$, where $a_i=a_{m-i}$ and $a_0\neq 0$. Now, if $m$ is odd, $d_n(1)=\sum_{i=0}^m a_i =\sum_{i=0}^{\frac{m-1}{2}} 2a_i = 0  \mod 2$. But writing out $\frac{1}{\det(R)}\det(\tau R t+R'\tau')$, $c_n(t)=\pm\frac{det(\theta t+\theta')}{\det(\theta+\theta')}$, so $c_n(1)=\pm 1$, but by our assumption, $\det (R) = 1 \mod 2$. This yields a contradiction.  So $m$ must be even. In this case $p(1)=a_{m/2}+\sum_{i=0}^{\frac{m}{2}-1} 2a_i$, and again because $c_n(1)=\pm 1$ and $\det (R)= 1 \mod 2$, $a_{m/2}$ must be odd. Consequently, $d_n(-1)=a_{m/2}+\sum_{\overset{i=0}{i \text{ even}}}^{\frac{m}{2}-1} 2a_i-\sum_{\overset{i=0}{i \text{ odd}}}^{\frac{m}{2}-1} 2a_i$ must also be odd. So, $c_n(-1)=\frac{d_n(-1)}{\det (R)}$ must be odd. But $c_n(-1)=\frac{1}{\det(R)}\det(-\tau R +R'\tau')$. The last determinant is that of a skew symmetric matrix and so must be of even dimension to be non-zero. This completes the argument.
\end{proof}

We next examine the relationship between $\theta$ and the Blanchfield pairing on the cokernel of $H_n(\td X)\to H_n(\td C)$ mod $\Z$-torsion.  Let us call this module $\bar H$ and  recall some facts from \cite[\S 3.6]{GBF1} (N.B. we have altered the notation from \cite{GBF1} in the hopes of introducing simpler and more consistent notation). It is shown there that for a disk knot  $L: D^{2n-1}\into D^{2n+1}$, $\bar H\otimes \Q\cong H_n(\td C;\Q)/(\ker p: H_n(\td C;\Q)\to H_n(\td C,\td X;\Q))$ is presented as a $\Gamma=\Q[\Z]$-module by the matrix $(-1)^{n+1}(R^{-1})'\tau R t-\tau'$ representing a map from $\bar E\otimes \Gamma \to \bar G\otimes \Gamma$. The only requirements assumed on the \emph{integral} bases of $\bar E$, $E$, $\bar G$, and $G$ are that those of $E$ and $\bar G$ are dual under the linking pairing $L'$ and similarly for $\bar E$ and  $G$ with $L'$. In this case, $R$ is simply the matrix of $p_*:\bar E\to \td E$. Also with respect to these integral bases (which induce an integral basis for $\bar H\otimes \Q$), the matrix of the self-Blanchfield pairing on $\bar H$ is given by  $\frac{t-1}{(R^{-1})'\tau -(-1)^{n+1}t\tau'R^{-1}}$. 

Let us demonstrate that the same matrix $M=(-1)^{n+1}(R^{-1})'\tau R t-\tau'$ in fact presents $\bar H$ as a $\Lambda$-module. 

\begin{proposition}
The  matrix $M=(-1)^{n+1}(R^{-1})'\tau R t-\tau'$ presents the $\Lambda$-module $\bar H$, which is the cokernel of  $H_n(\td X)\to H_n(\td C)$ modulo its $\Z$-torsion.
\end{proposition}
\begin{proof}
Let us denote $\mc{E}=\cok (H_n(F)\to H_n(V))$, $\mc{G}=\cok(H_n(S^{2n}-F)\to H_n(D^{2n+1}-V))$ and $\mc{H}=\cok(H_n(\td X)\to H_n(\td C))$ and consider the following commutative diagram:

\begin{diagram}
0 &\rTo  & H_n(F)\otimes \Lambda & \rTo & H_n(S^{2n}-F)\otimes \Lambda &\rTo & H_n(\td X)&\rTo & 0\\
   &           & \dTo                                   &          &                \dTo                                     &        &      \dTo    \\
0 &\rTo  & H_n(V)\otimes \Lambda & \rTo & H_n(D^{2n+1}-V)\otimes \Lambda &\rTo & H_n(\td C)&\rTo & 0\\                                          
   &           & \dOnto                                   &          &                \dOnto                                    &        &      \dOnto   \\
  &\rTo  & \mc{E}\otimes \Lambda  & \rTo^{\phi} & \mc{G}\otimes \Lambda                  &\rTo^{\eta} & \mc{H}    &\rTo & 0\\
  &           & \dOnto^{\psi}                                  &          &                \dOnto^{\xi}                                    &        &      \dOnto^{\zeta}    \\
  &           & \bar{E}\otimes \Lambda  & \rTo^f & \bar{G}\otimes \Lambda                  &\rTo^g & \bar{H}    \\
  &           & \dInto                                   &          &               \dInto                                  &        &      \dInto     \\
  &           & \bar{E}\otimes \Gamma  & \rTo & \bar{G}\otimes \Gamma                  &\rTo & \bar{H}\otimes_{\Z} \Q    &\rTo & 0 &.
\end{diagram}

The top row comes from the usual Mayer-Vietoris sequence for constructing an infinite cyclic cover of a knot by cutting and pasting along the Seifert surface. This sequence splits into short exact sequences using the fact that $H_n(\td X)$ is of type K; see \cite[p. 43]{L77}. The second row is also exact and arises from the same considerations applied to the disk knot. The third row is from the serpent lemma as this row consists of cokernels, and it is also exact; note that since $\Lambda$ is a free abelian group, tensoring by $\otimes_{\Z} \Lambda$ preserves exactness. The next row comes from killing all $\Z$-torsion, and the last row comes by taking the tensor product $\otimes_{\Z}\Q$. Note that $\Gamma=\Lambda\otimes_{\Z}\Q$. The bottom row is exact since it is the tensor product of the exact third row with $\Q$. The maps to the bottom row are injective since there is no $\Z$-torsion in the fourth row. Also the maps from the second row to the third and from the third to the fourth are clearly onto. Our goal is to show that the fourth row is short exact and presents $\bar H$ by the matrix $M$.

First, we verify that the map  $f: \bar{E}\otimes \Lambda  \to \bar{G}\otimes \Lambda$ is really our matrix $M$. To see this, we observe that the maps $d: H_n(F)\otimes \Lambda \to H_n(S^{2n}-F)$ and $d':  H_n(V)\otimes \Lambda \to H_n(D^{2n+1}-V)\otimes \Lambda$ come from the Mayer-Vietoris sequences and so both have the form $i_{-*}\otimes t-i_{+*}\otimes 1$; see \cite{L66} and \cite{GBF1}. By commutativity, this is then the form of the map $\bar{E}\otimes \Lambda  \to \bar{G}\otimes \Lambda$ under the quotients to the cokernels and mod torsion. But now $\{\bar \delta_j\otimes 1\}_{i=1}^m$  is a $\Lambda$-module basis for $\bar E\otimes \Lambda$, and we know by the definitions in Section \ref{S: basics} that
\begin{align*}
(i_{-*}\otimes t-i_{+*}\otimes 1)(\bar \delta_j\otimes 1) &= i_{-*}(\bar \delta_j)\otimes t-i_{+*}(\bar \delta_j)\otimes 1\\
&=\sum_i (\sigma_{ij}\bar \beta_i \otimes t - \lambda_{ij}\bar \beta_i \otimes 1).
\end{align*}
So we see that $f$ is represented by the matrix $\sigma\otimes t -\lambda\otimes 1$. But notice that we can now use the dualities discussed in Section \ref{S: basics} to see that, e.g. $\sigma_{kj}=L'(\alpha_k, i_{-*}\bar \delta_j)$.  So we can now apply the various properties of linking pairings as discussed in \cite{GBF1}. These properties hold integrally as well as rationally, and we can duplicate the arguments of \cite{GBF1} to see that $\sigma=\mu'$ and $\tau=\lambda'$. Also as in \cite{GBF1}, $\mu R=(-1)^{n+1}R'\tau'$. So we see that $f$ is indeed represented by   $(-1)^{n+1}(R^{-1})'\tau R \otimes t    -\tau'$, which we abbreviate as $M=(-1)^{n+1}(R^{-1})'\tau R  t    -\tau'$. This argument simply demonstrates that the rational presentation matrices obtained in \cite{GBF1} are really just this integral matrix tensored with $1\in \Q$. 

So now let us see that the fourth row of the diagram is exact: The map $f$ is injective because $\det(M(1))=\pm 1$, so $f$ is an injective $\Lambda$-module morphism. The map $g$ is onto by some easy diagram chasing. The composite $fg=0$ since the fourth row injects into the exact row below it. Finally, to see that $\ker g\subset \im f$, suppose that $0\neq x\in \ker g$. By the surjectivity of $\xi$, we know that $x=\xi(z)$ for some $z \in \mc{G}\otimes \Lambda  $. By commutativity, $\zeta\eta(z)=0$. This implies that $\eta(z)\in T(\mc{H})$, the $\Z$-torsion subgroup of $\mc{H}$. So there is an $m\in \Z$, $m\neq 0$ such that $m\eta(z)=\eta(mz)=0$, which implies that $mz\in \im(\phi)$, say $mz=\phi (y)$. So then $f\psi(y)=\xi\phi(y)=mx\neq 0$ because $\bar{G}\otimes \Lambda $ has no $\Z$-torsion. So $mx\in \im (f)$. But now consider the quotient $\Lambda$-module $A=(\bar{G}\otimes \Lambda) /\im(f)\cong (\bar{G}\otimes \Lambda)/M$. By the proof of \cite[Lemma 2.1]{T73}, $A$ is $\Z$-torsion free. So if $mx\in \im (f)$, we must also have $x\in \im(f)$. This completes the proof.
\end{proof}

In particular, this proposition implies that the basis $\bar \beta_i\otimes 1\in \bar G\otimes \Lambda$ spans $\bar H$. So the matrix $\frac{t-1}{(R^{-1})'\tau -(-1)^{n+1}t\tau'R^{-1}}$ also represents the \emph{integral} pairing 
$\bar H\times \bar H \to Q(\Lambda)/ \Lambda$; see \cite[\S 3.6.3]{GBF1} for the geometry that gives this formula, and observe that the calculation there is also made with respect to integral bases. Of course this matrix also represents the rational pairing $\bar H\otimes \Q\times \bar H\otimes \Q\to Q(\Gamma)/\Gamma$ that we obtain by tensoring everything  with $\Q$. Note, by the way, that $Q(\Lambda)=Q(\Gamma)$, both being the field of rational functions. 

Our goal now is to prove the following theorem:

\begin{theorem}\label{T: isometry implies cobordism}
Let $\theta_1$ and $\theta_2$ be Seifert matrices for disk knots $L_1,L_2: D^{2n-1}\into D^{2n+1}$. Suppose that $L_1$ and $L_2$ have isometric Blanchfield-self pairings on $\bar H_1 \cong \bar H_2$. Then $\theta_1$ and $\theta_2$ are integrally cobordant, in fact rationally S-equivalent.
\end{theorem}

Since Corollary \ref{C: rat gives int} tells us that it is enough to study Seifert matrices up to \emph{rational cobordism},  we can perform \emph{rational} changes of basis, maintaining the dualities with respect to the rational pairings $L'$ and $L''$, such that $R$ becomes the identity matrix. Then, rewriting the above matrices  using $\theta= R'\tau'$, we obtain the presentation matrix $(-1)^{n+1}\theta't-\theta$ for $\bar H \otimes \Q$ and rational Blanchfield pairing matrix $\frac{t-1}{\theta'-(-1)^{n+1}t\theta}$. These represent the same \emph{rational} module and pairing we started with up to isometry since we have only performed rational changes of basis on $\bar E\otimes \Q$ and $\bar G\otimes \Q$ and hence to $\bar H\otimes \Q$.  This nice new form puts us in position to use some slightly modified  machinery of Trotter \cite{T73}, though the transition from integral to rational will simplify things considerably. Note that for the duration of the proof of the theorem, we suspend our standard rule and allow the term \emph{Seifert matrix} to refer also to this new rational  $\theta$ obtained from the integral one by a rational change of bases.

We will need the notion of rational S-equivalence. For two square rational matrices $A$ and $B$, we say that $A$ is a rational row enlargement of $B$ and $B$ is a rational row reduction of $A$ if 
\begin{equation*}
A=\begin{pmatrix}
0&0&0\\
1&x& u\\
0&v&B
\end{pmatrix},
\end{equation*}
where $x$ and $1$ are rational numbers, $v$ is a column vector, and everything else is made to make the matrix square. Rational column enlargements and reductions are defined similarly with the transposed form. Rational S-equivalence is then the equivalence relation generated by rational row and column enlargements and reductions and by rational congruence.

\begin{lemma}\label{L: non-sing S-equiv}
For any  disk knot Seifert matrix $\theta$,   either  $\theta$ is rationally S-equivalent to a rationally nonsingular matrix or  $(-1)^{n+1}\theta't-\theta$ presents the $0$ $\Gamma$-module. 
\end{lemma}
\begin{proof}
It is shown on pages 484-485 of \cite{T62} that given an integral matrix $V$ with zero determinant and such that $\det(V-V')\neq 0$, then $V$ is integrally congruent to a matrix of the form
\begin{equation*}
\begin{pmatrix}
0&0&0\\
-1&0&0\\
0&q&W
\end{pmatrix},
\end{equation*}
where $W$ has dimensions $2$ less than those of $V$ and all other non-zero entries are integers. The same argument given there works, however, with $\det(V+V')\neq 0$ (the skew symmetry of $V-V'$ is mentioned but never used) and with ``integral'' replaced by ``rational'' at all steps. Note that, as usual, a rational matrix is considered rationally unimodular as long as its determinant is non-zero. This matrix demonstrates an S-equivalence between $V$ and $W$ (the $-1$ can be changed to a $1$ by a rational congruence).

So now, as in \cite[Lemma 1.4]{T73}, we can apply this process to $\theta$ inductively to reduce $\theta$ in dimension. Eventually we will obtain either  a nonsingular matrix or  a matrix of the form $(0)$ or $\begin{pmatrix}0&0\\1&x
\end{pmatrix}$. The first form is impossible since it follows from an elementary computation that if $V\pm V'$ is nonsingular then so is $W\pm W'$, and the second form presents the $0$ $\Gamma$-module when plugged into the formula. 
\end{proof}

\begin{remark}\label{R: nonsing}
We observe that if $V$ is nonsingular, then $\Delta=\det(tV\pm V')$ is a polynomial of degree equal to the dimension of $V$ and with non-zero constant term. The latter claim is clear by plugging in $t=0$. For the former, the nonsingularity implies that $\det(tV\pm V')=\det(V) \det(tI\pm V^{-1}V')$, which clearly has a term of the required degree.
\end{remark}

We will see in the next lemma that two rationally S-equivalent matrices present the same $\Gamma$-module. 

\begin{lemma}\label{L: S isometry}
If $\theta_1$ and $\theta_2$ are rationally S-equivalent, then they determine isometric $\Gamma$-modules with self-Blanchfield pairings.
\end{lemma}
\begin{proof}
The proofs of Lemmas 1.4 and 1.2 of \cite{T73} apply rationally. It should be  noted that our presentation matrix and pairing matrix defer slightly  from those in \cite{T73}. One reason is that we employ a different convention for turning a matrix into a pairing matrix (we use $a_1\times a_2\to a_1'M\bar a_2$, while Trotter uses $a_1\times a_2\to \bar a_2'Ma_1$. The other difference is the appearance of $\frac{1}{t-1}$ in Trotter's presentation matrices, but, as noted on \cite[p. 179]{T73}, these make no difference as multiplication by $t-1$ is an automorphism of knot modules. So the translation to Trotter's algebraic language from the topological  language can be made via some isomorphisms and convention switches, and so his results apply to our case. (One should also note carefully that what he calls $\Lambda$ is our $\Z[t,t^{-1},(1-t)^{-1}]$, while our $\Lambda$ is there denoted $\Lambda_0$.) 
\end{proof}

We will next need to consider Trotter's \emph{trace} function \cite{T73}: Since the rational functions, i.e. elements of  $Q(\Lambda)$, can be written in terms of partial fractions, $Q(\Lambda)$  splits over $\Q$ into the direct sum of  $\Gamma[(1-t)^{-1}]$ and the subspace $P$ consisting of $0$ and proper fractions with denominators prime to $t$ and $1-t$. The trace $\chi$ is then defined as the $\Q$-linear map to $\Q$ determined by $\chi(f)=f'(1)$ if $f\in P$ and $0$ if $f\in \Gamma[(1-t)^{-1}]$. The $'$ here denotes derivative with respect to $t$. This then induces a map $Q(\Gamma)/\Gamma\cong Q(\Lambda)/\Gamma\to \Q$. In particular, by composing $\chi$ with the Blanchfield pairing, one obtains a rational scalar form $\bar H\otimes \Q\times \bar H\otimes \Q\to \Q$. 

It is clear that two Seifert matrices that induce isometric Blanchfield forms induce isometric rational scalar forms.

Now by \cite[Lemma 2.7b]{T73}, for $f\in P$, $\chi((t-1)f)=f(1)$. And also, as in \cite[Lemma 2.10]{T73} and our Remark \ref{R: nonsing},  $\Delta$ has degree equal to the dimension of $\theta$ and non-zero constant term, plus we know it is prime to $(t-1)$, so by Cramer's rule, each term in $(\theta'-(-1)^{n+1}t\theta)^{-1} $ lies in $P$. Thus $\chi$ applied to $\frac{t-1}{\theta'-(-1)^{n+1}t\theta}$ is give by evaluation of $\frac{1}{\theta'-(-1)^{n+1}t\theta}$ at $1$, so we just get $\frac{1}{\theta'-(-1)^{n+1}\theta}$ as the matrix of the rational scalar pairing. 

It now follows as in the proof of  \cite[2.11]{T73}, using  \cite[2.5 and 2.10]{T73}, which also hold rationally, that a choice of basis in an  isometry class of a finitely generated $\Lambda[(t-1)^{-1}]$-module $ H_0$ with a rational scalar form \emph{determines} a ``Seifert matrix'' $\theta_0$ and that our given $\bar H\otimes \Q$ with rational scalar form is isometric to $ H_0$ if and only there is a basis for $\bar H\otimes \Q$ with respect to which its Seifert matrix $\theta$ is equal to $\theta_0$: The existence of an isometry implies that there are  bases with respect to which both scalar forms have the same matrix $S$ of \cite{T73}, and, with respect to these bases, $(1-t)^{-1}$ acts by the same matrix $\gamma$, but then the equations in \cite{T73} determine both $\theta_0$ and $\theta$ by $\gamma S^{-1}$. Finally, by \cite[Prop. 2.12]{T73}, this implies that two rationally nonsingular Seifert matrices determine isometric rational scalar forms if and only if they are rationally congruent. 

We can now complete the proof of Theorem \ref{T: isometry implies cobordism}. By hypothesis $\theta_1$ and $\theta_2$ determine isometric Blanchfield forms, hence they induce isometric scalar forms. Furthermore, by Lemma \ref{L: non-sing S-equiv}, $\theta_1$ and $\theta_2$ are rationally S-equivalent to Seifert forms, say $\hat \theta_1$ and $\hat \theta_2$, respectively, that are rationally nonsingular and which, by Lemma \ref{L: S isometry}, still determine isometric scalar forms. By the immediately preceding discussion, $\hat \theta_1$ and $\hat \theta_2$ are rationally congruent. It follows that $\theta_1$ and $\theta_2$ are rationally S-equivalent and hence, in particular, cobordant as seen in the proof of Proposition \ref{L: s-equivalence cobordism}. 

This completes the proof of Theorem \ref{T: isometry implies cobordism}. \qedsymbol 

\vskip1cm

The relationships we have just established between Seifert matrices and Blanchfield pairings turn out to be just what we need to realize rational cobordism classes of Seifert matrices.

\begin{theorem}\label{T: real simple}
Let $\theta$ be any square matrix satisfying the necessary conditions to be the integral Seifert matrix of a disk knot $D^{2n-1}\into D^{2n+1}$, i.e. such that 
\begin{enumerate}
\item $R_{\theta}=-\theta'+(-1)^{n+1}\theta$ is  nondegenerate, and

\item $\tau_{\theta}=\theta'(-\theta'+(-1)^{n+1}\theta)^{-1}$ and $\mu_{\theta}=(-\theta+(-1)^{n+1}\theta')^{-1}\theta'$ are integral matrices.
\end{enumerate}
Then for any $n>2$, there is a disk knot $D^{2n-1}\into D^{2n+1}$ whose Seifert matrix is cobordant to $\theta$.
\end{theorem}
\begin{proof}
Given such a $\theta$, it determines a $\Lambda$-module $\bar H$ with a $(-1)^{n+1}$-Hermitian pairing to $Q(\Lambda)/\Lambda$ by the matrices  $(-1)^{n+1}(R^{-1})'\tau R t-\tau'$ and $\frac{t-1}{(R^{-1})'\tau -(-1)^{n+1}t\tau'R^{-1}}$ as in the discussion earlier in this section (see also \cite[\S 3.6.3]{GBF1}). Note that $\bar H$ is $\Z$-torsion free by the same arguments as in \cite[Lemma 2.1]{T73}. By \cite[Proposition 3.21]{GBF1}, there exists a simple disk knot $L$ realizing this module and pairing with $\bar H=H_n(\td C)$ and also with simple boundary knot such that $H_{n-1}(\td X)$ is $\Z$-torsion. By Theorem \ref{T: isometry implies cobordism}, any Seifert matrix for $L$ is cobordant to our given $\theta$; in fact it is rationally S-equivalent to it. 
\end{proof}

So, at this point we have demonstrated that, for $n>2$, every cobordism class can be realized by 1) showing that a potential  Seifert matrix determines a Blanchfield pairing, 2) constructing every possible Blanchfield pairing, and 3) showing that a Blanchfield pairings determine its Seifert matrices up to rational S-equivalence. So by constructing every possible pairing, we construct every possible cobordism class. However, we have not said anything yet about what boundary knots we get. The constructions of Theorem \ref{T: real simple} give only simple disk knots whose boundaries are simple sphere knots and such that $H_{n-1}(\td X)$ is $\Z$-torsion  (this follows from the construction in \cite[Prop. 3.21]{GBF1} and the construction in \cite[\S 12]{L77} that it is modeled after). Such sphere knots are called \emph{finite simple}. In this special case, we can say a lot immediately. We will show in Section \ref{S: blanch to F-L} below that in this situation  the Blanchfield pairing on $H_n(\td C)$ completely determines the Farber-Levine torsion pairing on $H_{n-1}(\td X)$. In fact, we will prove the following theorem:

\begin{theorem}[Corollary \ref{C: Seifert gives FL}]\label{T: S to FL}
For a simple disk knot $L:D^{2n-1}\subset D^{2n+1}$, the $\Lambda$-module $T_{n-1}(\td X)$ and its Farber-Levine $\Z$-torsion pairing are determined up to isometry by any Seifert matrix for $L$.
\end{theorem}

In this situation, we will say that \emph{the Seifert matrix induces the Farber-Levine pairing}.

 We can now apply the following theorem of Kojima \cite{Koj} (which we have translated into our language):

\begin{theorem}[Kojima]\label{T: Kojima}
Suppose that  $K_0$ and $K_1$ are two finite simple sphere knots $S^{2n-2}\to S^{2n}$, $n\geq 5$,  $H_{n-1}(\td X_0)\cong H_{n-1}(\td X_1)$  contains no $2$-torsion, and the Farber-Levine pairings on $H_{n-1}(\td X_0)$ and $H_{n-1}(\td X_1)$ are isometric, then $K_0$ and $K_1$ are isotopic knots.
\end{theorem}

Putting this theorem together with the results of Section \ref{S: blanch to F-L}, quoted above, we see that, for $n\geq 5$, the following statement holds: if a Blanchfield pairing on $H_n(\td C)$ induces a $T_{n-1}(\td X)$ with no $2$-torsion, then this Blanchfield pairing determines a unique finite simple sphere knot $S^{2n-2}\into S^{2n}$ which must be the boundary knot of any simple disk knot  possessing this Blanchfield pairing and having a  finite simple boundary knot.
In particular then, since Seifert matrices determine Blanchfield pairings, the Seifert matrix of a simple disk knot with finite simple boundary knot determines the boundary knot uniquely, so long as $H_{n-1}(\td X)=T_{n-1}(\td X)$ has no $2$-torsion. 

We can now immediately generalize this to prove the following theorem about realizability of cobordism classes of Seifert matrices for more arbitrary boundary knots:

\begin{theorem}\label{T: realize}
Let $K:S^{2n-2}\into S^{2n}$, $n\geq 5$, be a sphere knot with complement $X$ such that $T_{n-1}(\td X)$ contains no $2$-torsion. Then there exists a disk knot $L:D^{2n-1}\into D^{2n+1}$ with boundary knot $K$ and with Seifert matrix in a given cobordism class $[\theta]$ if and only if there is an integral matrix $\theta$ in the class such that 
\begin{enumerate}
\item $R_{\theta}=-\theta'+(-1)^{n+1}\theta$ is  nondegenerate,

\item $\tau_{\theta}=\theta'(-\theta'+(-1)^{n+1}\theta)^{-1}$ and $\mu_{\theta}=(-\theta+(-1)^{n+1}\theta')^{-1}\theta'$ are integral matrices, and
\item the Farber-Levine pairing induced by $\theta$ is isometric to the Farber-Levine pairing on $T_{n-1}(\td X)$. 
\end{enumerate}
\end{theorem}
\begin{proof}
Suppose we have such a knot $L$ and its cobordism class of Seifert matrices $[\theta]$. We show that there is a Seifert matrix in the cobordism class satisfying the listed properties:
We know that the first two requirements are always necessary for a Seifert matrix. For the third, recall that by Theorem \ref{T: cob to simple}, any disk knot is cobordant rel boundary to a simple simple disk knot, and by Theorem \ref{T: main cobordism}, any two such disk knots have cobordant Seifert matrices. By Theorem \ref{T: S to FL}, any Seifert matrix of a simple disk knot determines the Farber-Levine pairing on $T_{n-1}(\td X)$ of the boundary knot up to isometry. So there is a Seifert matrix in the cobordism class $[\theta]$ that induces the correct Farber-Levine pairing (up to isometry).

Conversely, given a $\theta$ that meets the above requirements, Theorem \ref{T: real simple} and its proof  assure us that we can construct a simple disk knot $L_1$ with finite simple boundary whose Seifert matrices fall in the cobordism class $[\theta]$ of $\theta$ and induce the given Farber-Levine pairing on the boundary knot. Now let $L_0$ be any simple disk knot with our given $K$ as boundary. Such a knot always exists since $K$ is null-cobordant by its dimensions and \cite{Ke}, and there is  a cobordism rel boundary of any disk knot to a simple disk knot by Theorem \ref{T: cob to simple}. Let $\theta_0$ be any Seifert matrix of $L_0$, and note that $\theta_0$ determines the Farber-Levine pairing on $T_{n-1}(\td X)$. Also, again by Theorem \ref{T: real simple}, there is a simple disk knot with torsion simple boundary $\theta_0$ whose Seifert matrices fall in the  cobordism class $[\theta_0]$ and induce the given Farber-Levine pairing. Since $L_1$ and $L_0$ are both simple disk knots with torsion simple boundaries $K_1$ and $K_0$ and since the boundary modules $H_{n-1}(\td X_0)$ and $H_{n-1}(\td X_1)$ are Farber-Levine isometric by construction and contain no $2$-torsion by assumption, Kojima's Theorem \cite{Koj} implies that $K_0$ and $K_1$ are isometric. So now let us form the sphere knot $\mc{K}=L_1\cup_{K_0}-L_0$. By Theorem \ref{T: combo knot}, the Seifert matrix of $\mc{K}$ is  cobordant to $\theta\boxplus -\theta_0$.  Finally, we form the connected sum away from the boundary $L=L_0\#\mc{K}$. Then $L$ has Seifert matrix  cobordant to $\theta_0\boxplus (\theta\boxplus -\theta_0)=\theta$, and it is our desired knot. 
\end{proof}

We note that the statement of the theorem only guarantees that some element in the cobordant class determines the proper Farber-Levine pairing, not all elements. This is really the best that can be hoped for since given an arbitrary disk knot, it is possible that $T_{n-1}(\td X)$ may not be in the image of $\bd_*$ or there may be elements in $T_{n-1}(\td X)$  that are in the image of $T_{n}(\td C,\td X)$. The Farber-Levine pairing on such elements clearly won't be determined by the Seifert matrix. However, as noted in the proof, there is always a cobordism rel boundary to a simple disk knot for which the entirety of the Farber-Levine pairing \emph{is} determined by the Seifert matrix, and we know that such a cobordism keeps the Seifert matrix in its cobordism class. While this argument shows that a  cobordism class does not determine a Farber-Levine pairing, we make the following conjecture:
\paragraph{Conjecture}
The  cobordism class of any integer matrix satisfying 
\begin{enumerate}\item $R_{\theta}=-\theta'+(-1)^{n+1}\theta$ is  nondegenerate,

\item $\tau_{\theta}=\theta'(-\theta'+(-1)^{n+1}\theta)^{-1}$ and $\mu_{\theta}=(-\theta+(-1)^{n+1}\theta')^{-1}\theta'$ are integral matrices
\end{enumerate}
determines a unique element in the Witt group  of $\Z$-linear conjugate self-adjoint $(-1)^{n+1}$-symmetric nonsingular pairings to $\Q/\Z$ on finite $\Lambda[(t-1)^{-1}]$-modules.  

Our realization theorem makes no conclusions about knots for which $T_{n-1}(\td X)$ possesses $2$-torsion. This is because finite simple even-dimensional sphere knots are not determined entirely by their Farber-Levine pairings, and so the previous proof breaks down; we can not apply the theorem of Kojima. It was shown by Farber in a series of papers culminating in \cite{Fa84a, Fa84b} (see also \cite{Fa83}) that in this case there is also an even-torsion pairing on the stable homotopy groups $\sigma_{n+1}(\td X)$ that plays a role in the classification. In fact, Farber shows that such knots are classified completely by the algebraic invariants in their \emph{$\Lambda$-quintets}. It remains unclear whether the Seifert matrices and/or Blanchfield pairings of a simple disk knot are sufficient to determine the $\Lambda$-quintets of their boundary knots, so we can not yet broaden Theorem \ref{T: realize} to include realizability for all knots. An alternative procedure would be to show that all knots constructed in Theorem \ref{T: real simple} that  give the same Farber-Levine pairing on the boundary just happen to  have the same actual boundary knot. If so, the proof of Theorem \ref{T: realize} would apply without the need to invoke a broader classification theorem. However, we have not yet been able to establish this either.

\section{Blanchfield pairings determine Farber-Levine pairings}\label{S: blanch to F-L}

In this section, we will establish that for a simple disk knots of odd dimension $D^{2n-1}\subset D^{2n+1}$,  the Farber-Levine $\Z$-torsion self-pairing $T_{n-1}(\td X)\otimes T_{n-1}(\td X)\to \Q/\Z$ is determined completely by the module $\bar H=\cok(H_n(\td X)\to H_n(\td X))$ and its self-Blanchfield pairing. This result is used in the previous section in conjunction with the main theorem of \cite{Koj} to recognize the boundary knots of knots we have constructed. 

We will begin by demonstrating that the module $H_{n-1}(\td X)$  and the Farber-Levine pairing on its submodule $T_{n-1}(\td X)$ are determined by the self-Blanchfield pairing on $H_n(\td C)$. This will be done initially by developing a formula relating the two pairings based upon the geometry of chains. Once this connecting formula is established, we will abstract to the purely algebraic situation and redefine the Farber-Levine pairing by a completely algebraic construction given $H_n(\td C)$ and its Blanchfield pairing. This will allow us to prove that the isometry class of the latter completely determines  the isometry class of the former. We then show that, in fact, $\bar H$, which algebraically corresponds to the quotient of $H_n(\td C)$ by its annihilating submodule, is sufficient to determine $T_{n-1}(\td X)$ and its Farber-Levine pairing.

To simplify things marginally,  observe that $\bd \td C =\td X \cup_{S^{2n-2}\times \R} D^{2n-1}\times \R$ so that, for $n\geq 2$,  the map induced by inclusions $H_{n-1}(\td X)\to H_{n-1}(\bd \td C)$ is an isomorphism  and $H_n(\td X)\to H_n(\td C)$ is an epimorphism. It therefore follows from the five lemma applied to the exact sequences of the pairs that $H_{n}(\td C, \td X)\to H_{n}(\td C, \bd \td C)$ is an isomorphism. For $n=1$, $X\sim_{h.e.} S^1$, so $\td X \sim_{h.e.}*$. In this case there is no Farber-Levine pairing of interest, so we will shall always assume $n\geq 2$. We will work with $\bd \td C$ or $\td X$  as  convenient, but using these isomorphisms, we can assume that all relevant chains are actually contained in $\td X$.

For a simple disk knot $D^{2n-1}\into D^{2n+1}$,  $H_i(\td C)=0$ for $0<i<n$ due to the connectivity assumptions.  Now, as observed in \cite{L77} (and holding for  any regular covering of a compact piecewise-linear $n$-manifold with boundary), $H_{*}(\td C)\cong \overline {H_e^{2n+1-*}(\td C, \bd \td C)}$, the conjugate of the cohomology of the cochain complex $\Hom_{\Lambda}(C_*(\td C, \bd \td C),\Lambda)$. Similarly, $H_{*}(\td C, \bd \td C)\cong \overline{ H_e^{2n+1-*}(\td C)}$, the conjugate of the cohomology of the cochain complex $\Hom_{\Lambda}(C_*(\td C),\Lambda)$. It now follows from Proposition 2.4 of \cite{L77} and this generalization of Poincar\'e duality that there exist short exact sequences 
\begin{equation*}
\begin{CD}
0&@>>>& \Ext^2_{\Lambda}(H_{n-2}(\td C), \Lambda) &@>>> & \overline{H_{n+1}(\td C,\bd \td C)}&@>>> \Ext^1_{\Lambda}(H_{n-1}(\td C), \Lambda) &@>>> 0\\
0&@>>>& \Ext^2_{\Lambda}(H_{n-1}(\td C), \Lambda) &@>>> &\overline{ H_n(\td C, \bd \td C)}&@>>> \Ext^1_{\Lambda}(H_{n}(\td C), \Lambda) &@>>> 0.\\
\end{CD}
\end{equation*}
By the connectivity assumptions on $\td C$, these imply that $H_{n+1}(\td C,\bd \td C)=H_{n+1}(\td C,\td X)=0$, and 
$\overline{ H_n(\td C, \bd \td C)}\cong\Ext^1_{\Lambda}(H_{n}(\td C), \Lambda)$. Since  $H_n(\td C)$ is of type $K$ (it is finitely generated and $t-1$ acts as an automorphism), $\Ext^1_{\Lambda}(H_{n}(\td C), \Lambda)$ is $\Z$-torsion free by \cite[Prop. 3.2]{L77}, hence so is $H_n(\td C, \bd \td C)\cong H_n(\td C,\td X)$.

So we have shown that there is an exact sequence of Alexander modules 
\begin{equation*}
\begin{CD}
0&@>>>&H_{n}(\td X)& @>>>& H_{n}(\td C)&@>p_*>>& H_{n}(\td C, \td X) &@>\bd_*>> & H_{n-1}(\td X) &@>>>& 0
\end{CD}
\end{equation*}
and that $H_{n}(\td C, \td X)$ has no $\Z$-torsion.
We seek first to determine how the self-Blanchfield pairing on $H_n(\td C)$ determines the Farber-Levine $\Z$-torsion pairing $[\,,\,]: T_{n-1}(\td X)\otimes T_{n-1}(\td X)\to \Q/\Z$, where $T_{n-1}(\bd X)$ is the $\Z$-torsion subgroup of $H_{n-1}(\td X)$.

We begin by recalling the constructions of the various pairings involved. The following discussion integrates the relevant work from papers of Blanchfield \cite{B57} and Levine \cite{L77} and adapts it, where necessary, to the case of disk knots.   

We can assume that $\td C$, the infinite cyclic cover of the exterior of the disk knot $L: D^{N-2}\subset D^N$, is triangulated equivariantly so that $C_*(\td C, \bd \td C)$ is a free left $\Lambda$-module with basis given by the cells of $C$ not in $\bd C$. Then $C_*(\td C)$ can be taken as the free left $\Lambda$-module with basis given by the dual cells to the given triangulation of $C$ \cite{L77}. One then defines an intersection pairing of left $\Lambda$-modules to $\Lambda$ at the chain level by setting $a\cdot b=\sum_i S(a,t^ib)t^i$ for $a\in C_i(\td C)$, $b\in C_{N-i}(\td C, \bd \td C)$, where $S$ is the ordinary intersection pairing of chains. If we use an overline $\bar {}$ to denote the antiautomorphism on $\Lambda$ determined by $\bar t=t^{-1}$,  this pairing satisfies the following properties  \cite{B57}:
\begin{enumerate}
\item $(x+y)\cdot z=x\cdot z+y\cdot z$
\item $x\cdot( y+z)=x\cdot y+x\cdot z$
\item $(\gamma x)\cdot  (\bar \delta y)=\gamma\delta (x\cdot y)$ for $\gamma, \delta\in \Lambda$
\item $x\cdot  \bd y=(-1)^i\bd x \cdot y$ for $x\in C_i(\td C)$, $y\in C_{N+1-i}(\td C, \bd \td C)$
\item There exist dual bases $\{x_j\}\subset C_i(\td C)$ and $\{y_k\}\subset C_{N-i}(\td C, \bd \td C)$ such that $x_j\cdot y_k=\delta_{jk}$, the Kronecker delta.  
\end{enumerate}
These properties ensure that the pairing  $\cdot $ descends to a well-defined pairing of homology modules. It also follows from the properties of the ordinary intersection form on a manifold that if $x\in H_i(\td C)$ and $y\in H_{N-i}(\td C)$, then $x\cdot p_*(y)=(-1)^{i(N-i)}\overline{y\cdot p_*(x)}$, where $p_*: H_n(\td C) \to H_n(\td C,  \bd \td C)$.

From here, it is possible to define a linking pairing (the Blanchfield pairing) $V: W_i(\td C)\otimes W_{N-1-i}(\td C, \bd \td C)\to Q(\Lambda)/\Lambda$, where $W_i(X)$ is the submodule of weak boundaries of $C_i(X)$ and $Q(\Lambda)$ is the field of rational functions. If $a\in W_i(\td C)$, $b\in W_{N-1-i}(\td C, \bd \td C)$ and $A\in C_{i+1}(\td C)$ with $\bd A=\alpha a$ for some $\alpha\in \Lambda$, then $V(a,b)=\frac{1}{\alpha}A\cdot b$ by definition. Note that this linking number is well-defined to $Q(\Lambda)$ at the chain level. However, in order to descend to a well-defined map on homology classes with torsion, it is necessary to consider the image of $V$ in $Q(\Lambda)/\Lambda$. In the case of interest to us, the relevant pairing will be $V:H_n(\td C)\otimes H_{n}(\td C, \bd \td C)\to Q(\Lambda)/\Lambda$ when $N=2n+1$ (recall that both modules are $\Lambda$-torsion so all cycles weakly bound). By \cite[\S 5]{L77}, since $H_n(\td C,\bd \td C)$ is $\Z$-torsion free, the pairing is nonsingular in the sense that its adjoint provides an isomorphism $\overline{ H_{n}(\td C, \bd \td C)}\to \Hom_{\Lambda}(H_n(\td C), Q(\Lambda)/\Lambda)$ (the overline on $\bar H_{n}(\td C, \bd \td C)$ indicates that we take the module with the conjugate action of $\Lambda$ under the standard antiautomorphism, reflecting the fact that $V$ will be conjugate linear, since $\cdot$ is). This pairing determines a self-pairing $\langle \,,\,\rangle$ on $H_n(\td C)$ by $\langle a,b\rangle=V(a,p_*(b))$.  This pairing is $(-1)^{n+1}$-Hermitian, i.e. $\langle a,b\rangle=(-1)^{n+1}\overline{\langle b,a\rangle}$, and it is  nondegenerate on $\text{coim}(p_*)$. 

Now, it requires more work to define the Farber-Levine $\Z$-torsion pairing. Generally, these are pairings  $[\,,\,]: T_i(\td X)\otimes T_{N-i-2}(\td X)\to \Q/\Z$, where $T_j(\td X)$ is the $\Z$-torsion submodule of $H_j(\td X)$ and  $X$ has dimension $N$. We will specialize immediately to our case of interest $N=2n$, $i=n-1$. Note that this is a pairing on the torsion Alexander module of a locally-flat sphere knot, so we simply repeat Levine's construction from \cite{L77}. In fact, Levine begins with a sophisticated definition via homological algebra and then produces an equivalent geometric formulation. We will be more concerned with the geometric formulation, but there is one intermediate algebraic construction that remains necessary. We first need to choose two integers, but the final outcome will be independent of the choice modulo the restrictions on choosing. Let $m$ be a positive integer such that $mT_e^{n+1}(\td X)=0$, where $T_e^{n+1}(\td X)$ is the torsion subgroup of $H_e^{n+1}(\td X)$. By generalized Poincar\'e duality, $T_e^{n+1}(\td X)\cong \overline{T_{n-1}(\td X)}$, so $m$ kills this module as well. In fact, such an $m$ exists since $T_{n-1}(\td X)$ is finite by \cite[Lemma 3.1]{L77}.  Next, let $\Lambda_m=\Lambda/m\Lambda=\Z_m[\Z]$, and  let $\theta=\Lambda/(t^k-1)$, where $k$ is a positive integer chosen large enough so that $t^k-1$ annihilates $H_e^{n}(\td X;\Lambda_m)$. Such a $k$ exists since $H_e^{n}(\td X;\Lambda_m)\cong \overline{H_n(\td X;\Lambda_m)}$ by generalized Poincar\'e duality, and this module is also finite, again by \cite[Lemma 3.1]{L77} and the argument on the bottom of page 18 of \cite{L77}.  
Since $H_e^n(\td X;\Lambda_m)$ is finite and $t$ acts isomorphically, $t^k=1$ for some integer $k>0$. Hence $t^k-1$ annihilates the module for this choice of $k$. Note that $t^k-1$ also kills $H_n(\td X;\Lambda_m)$ since $\overline{t^k-1}=t^{-k}-1=-t^{-k}(t^k-1)$ and $t$ acts automorphically. By the same arguments, we can find  a $k$ such that $t^{k}-1$ annihilates $H_n(\td C;\Lambda_m)$. Since $H_{n+1}(\td C,\td X)=0$, we also get $H_n(\td C,\td X;\Lambda_m)=0$ because, as an abelian group, $H_{n+1}(\td C,\td X;\Lambda_m)=H_{n+1}(\td C,\td X)\otimes_{\Z} \Z_m$ (recall that $H_{n}(\td C,\td X)$ is $\Z$-torsion free).  So $H_n(\td X;\Lambda_m)$ maps monomorphically into $H_n(\td C;\Lambda_m)$ in the long exact sequence of the pair $(\td C,\td X)$ with $\Lambda_m$ coefficients, so this $k$ suffices to kill $H_{n}(\td X;\Lambda_m)$ as well. In other words, for any $k$ such that $t^k-1$ kills $H_n(\td C;\Lambda_m)$, the same choice of $k$ gives a $t^k-1$ that also kills $H_n(\td X;\Lambda_m)$.

The geometric part of the construction now finds a pairing $\{\,,\,\}: T_{n-1}(\td X)\otimes T_{n-1}(\td X)\to I(\theta)/\theta$, where $I(\theta)$ is the $\Lambda$-injective envelope of $\theta$. But for a finite $\Lambda$-module $A$, $\Hom_{\Lambda}(A,I(\theta)/\theta)= \Hom_{\Lambda}(A, (\Q\otimes_{\Z} \theta)/\theta)\overset{e}{\cong} \Hom_{\Z}(A, \Q/\Z)$, so it is possible to define $[\,,\,]$ as the composition of $\{\,,\,\}$ with these isomorphisms. 

The pairing $\{\,,\,\}$ can be described in the following geometric manner: Suppose that $z$ and $w$ are cycles representing elements of $T_{n-1}(\td X)$. By the choice of $m$,  $mz$ is null-homologous, so $mz=\bd z'$ for some $z'\in C_{n}(\td X)$. Then $(t^k-1)z'$ is null-homologous mod $m$ since we know that $t^k-1$ annihilates $H_n(\td X;\Lambda_m)$. Thus, we can write $(t^k-1)z'=\bd z'' +mz_0$ for some $z''\in C_{n+1}(\td X)$ and $z_0\in C_{n}(\td X)$. Then one sets $\{z,w\}$ to be the image of $(-z''\cdot w)/m$, which is in $\Gamma$, under the composition $\Gamma\to\Q\otimes \theta\subset I(\theta)\to I(\theta)/\theta$. It turns out that this pairing is independent of the choices involved and descends to a well-defined map on the homology torsion subgroups. See \cite{L77} for more details.

Using this geometric definition, we next show how the middle dimensional pairing $\{\,,\,\}: T_{n-1}(\td X)\otimes T_{n-1}(\td X)\to I(\theta)/\theta$ can be expressed in terms of the linking pairing $V: H_n(\td C)\otimes H_n(\td C,\bd \td C)$. In the following computations, all pairings are defined at the chain level, so there is no ambiguity. Consider cycles $z,w$ representing elements in $T_{n-1}(\td X)$. Let $z', z'', z_0 \in C_*(\td X)$ be as defined above. We need to reformulate $(-z''\cdot w)/m$, where the intersection product is that in $\td X$. 

 Since $\bd_*:H_{n}(\td C, \td X)\to H_{n-1}(\td X)$ is surjective, there exist chains  $X, Y\in C_n(\td C)$ such that $\bd X=z$ and $\bd Y=w$.
Then, the intersection number $z''\cdot w$ in $\td X$ is equal to the intersection number of $z''$ and $Y$ in $\td C$. This follows just as in the more standard case of intersection numbers for manifolds with boundary. 

Next, observe that $0=\bd^2 z''=\bd ((t^k-1)z'-m z_0)=(t^k-1)mz-m\bd z_0$, which implies that $\bd z_0=(t^k-1)z$. This also implies the important fact that $t^k-1$ annihilates $T_{n-1}(\td X)$ since $z$ is an arbitrary element of it. 
 Let $S=(t^k-1)X-z_0\in C_n(\td C)$. The chain $S$ is a cycle and so represents an element of $H_{n}(\td C)$. Since $H_n(\td C)$ is a finitely generate $\Lambda$-torsion module, there exists an element $\Delta\in \Lambda$ such that $\Delta H_n(\td C)=0$. So there exists a chain $R\in C_{n+1}(\td C)$ such that $\bd R=\Delta S$. Similarly, define the $n$-cycle $B=mX-z'$, and choose an $n+1$ chain $A$ in $\td C$  such that $\bd A=\Delta B$. 

Now $\bd (mR-(t^{k}-1)A)=m\bd R-(t^k-1)\bd A=m\Delta S-(t^k-1)\Delta B=\Delta(m(t^k-1)X-mz_0-(t^k-1)( mX-z'))=\Delta \bd z''$. Using the properties of  intersection forms we can see that $\Delta z''\cdot Y=( m R-(t^k-1) A)\cdot Y$. In fact,
$(\Delta z''-mR+(t^k-1)A)$ is a cycle in $C_{n+1}(\td C)$ and so represents a  homology class. Thus $(\Delta z''-mR+(t^k-1)A)\cdot Y$ is a well-defined element of $\Lambda$ under the intersection pairing $H_{n+1}(\td C)\otimes H_n(\td C,\td X)\to \Lambda$. But we know this pairing is $\Lambda$-linear in $H_{n+1}(\td C)$, and $H_{n+1}(\td C)$ is $\Lambda$-torsion. So this intersection must be $0\in \Lambda$.

Thus, since the intersection $z''\cdot w$ in $\td X$ is equal to the intersection number $z''\cdot Y$ in $\td C$, we compute
\begin{align*}
\frac{z''\cdot_{\td X} w}{m}=\frac{z''\cdot_{\td C} Y}{m}&=\frac{( m R-(t^k-1) A)\cdot Y}{m\Delta}\\
&=\frac{R\cdot Y}{\Delta}-\frac{t^k-1}{m}\frac{A\cdot Y}{\Delta}\\
&=V(S,Y)-\frac{t^k-1}{m}V(B,Y),
\end{align*}
and this establishes a formula for $\{z,w\}=\frac{-z''\cdot_{\td X}w}{m}$ in terms of the linking pairing $V$ under the projection to $I(\theta)/\theta$. Note that this formula is well-defined on passage to homology, since we know in this  case that $V$ is well-defined up to elements of $\Lambda$. So the first term of this expression is well-defined up to an element of $\Lambda$ and the second term up to elements of the form $\frac{t^k-1}{m}\lambda$, $\lambda\in \Lambda$. But all such elements are in the kernel of the composition $\Gamma\to\Q\otimes \theta\subset I(\theta)\to I(\theta)/\theta$. Note, however, that we are not free to conclude that the term $\frac{t^k-1}{m}V(B,Y)$ lies in this kernel. 

Since this construction yields the well-defined element $[z,w]$, it must be independent of the choices made in the construction, but we will also verify this below in the process of abstracting this pairing to a purely algebraic construction. 
Thus we will see that the torsion pairing $[\,,\,]$ on $T_n(\td X)$ is completely determined by the isometry class of the Blanchfield self-pairing on $H_n(\td C)$. So now we forget the geometry and abstract to a purely algebraic setting (though we will, perhaps confusingly, keep the geometric notation). Suppose we are given a $\Lambda$-torsion module of type K and possessing a   $(-1)^{n+1}$-Hermitian self-pairing $\langle \,,\,\rangle$. We will suggestively call this module $H_n(\td C)$. This determines a map $p_*: H_n(\td C)\to \overline{\Hom(H_n(\td C), Q(\Lambda)/\Lambda)}$, and we suggestively call the codomain module $H_n(\td C, \bd \td C)$. It is $\Z$-torsion free by \cite[Props. 3.2, 4.1]{L77}.  The map is determined by the formula $\langle x,y\rangle=V(x,p_*(y))$, where here $V$ is the canonical pairing $H_n(\td C)\times \overline{\Hom(H_n(\td C), Q(\Lambda)/\Lambda)}\to Q(\Lambda)/\Lambda$. This in turn yields the quotient module $\cok(p_*)$, which we denote $H_{n-1}(\td X)$. The $\Z$-torsion submodule $T_{n-1}(\td X)$ will be finite by \cite[Lemma 3.1]{L77}, and so we can choose an integer $m>0$ such that $m$ annihilates it. Similarly, we find a $k$ such that $t^k-1$ annihilates $H_n(\td C)\otimes_{\Z}\Z_m$ and $T_{n-1}(\td X)$. To see that such a $k$ exists, we make the expedient observation that since $H_n(\td C)$ is a $\Lambda$-module of type $K$ by assumption, there exists a sphere knot with $H_n(\td C)$ as one of its Alexander modules and $0$ as its succeeding $\Lambda$-module, and then the existence of an annihilating $t^k-1$ for the corresponding $H_n(\td C;\Lambda_m)\cong H_n(\td C)\otimes_{\Z}\Z_m$ follows as in the original arguments in \cite{L77}. 
Furthermore, $\overline{\Hom(H_n(\td C), Q(\Lambda)/\Lambda)}$ will be its dual Alexander module, hence also of type $K$, and so the quotient $H_{n-1}(\td X)$ is of type K. So $T_{n-1}(\td X)$ is finite and also subject to annihilation by $t^k-1$ for some $k$. Hence there exists some $k$ such that $t^k=1$ on both these modules.

We can fix $m$ and $k$ by choosing the smallest positive integers that satisfy these properties.  Then given $x,y\in T_{n-1}(\td X)$, we can \emph{define} $\{x, y\}\in I(\theta)/\theta$ by choosing $X, Y\in H_n(\td C,\td X)$ such that $\bd_*(X)=x$ and $\bd_*(Y)=y$ and $B,S\in H_n(\td C)$ such that $p_*(B)=mX$ and $p_*(S)=(t^k-1)X$, which is possible since $\bd_*(mX)=\bd_*((t^k-1)X)=0$. Then to get $\{x,y\}$ we simply take $-V(S,Y)+\frac{t^k-1}{m}V(B,Y)\in Q(\Lambda)/\Lambda$, where $V: H_n(\td C)\otimes H_n(\td C,\bd \td C)\cong H_n(\td C) \otimes \overline{\Hom (H_n(\td C), Q(\Lambda)/\Lambda)}\to Q(\Lambda)/\Lambda$ is the natural Hermitian evaluation, and compose with the map $Q(\Lambda)/\Lambda\to I(\theta)/\theta$ induced by the commutative diagram 
\begin{equation*}
\begin{CD}
\Lambda &@>>> & \theta\\
@VVV&&@VVV\\
Q(\Lambda)&@>>>& I(\theta)
\end{CD}
\end{equation*}
whose vertical maps are inclusions (see \cite{L77}). To get $[x,y]$, we then follow \cite{L77} and apply the isomorphisms $\Hom_{\Lambda}(A,I(\theta)/\theta)= \Hom_{\Lambda}(A, (\Q\otimes_{\Z} \theta)/\theta)\overset{e}{\cong} \Hom_{\Z}(A, \Q/\Z)$ for a finite $\Lambda$-module $A$.

Let us show algebraically that this process is independent of choices of $X$, $Y$, $B$, and $S$. Continuing to mirror the geometric language, we let $\bd_*$ denote the quotient $H_{n}(\td C,\bd \td C)\to H_{n-1}(\td X)=\cok(p_*)$. Suppose that instead of $Y\in H_n(\td C,\bd \td C)$, we choose $Y'$ so that $\bd_* Y=\bd_*Y'=w$. Then $\bd_*(Y-Y')=0$ so $Y-Y'=p_*(U)$ for some  $U$ in $H_n(\td C)$. Then
\begin{align*}
V(S,Y)&-\frac{t^k-1}{m}V(B,Y)-(V(S,Y')-\frac{t^k-1}{m}V(B,Y'))\\
&=V(S,Y-Y')-\frac{t^k-1}{m}V(B,Y-Y')\\
&=V(S,p_*(U))-\frac{t^k-1}{m}V(B, p_*(U))\\
&=\langle S,U\rangle-\frac{t^k-1}{m}\langle B,U\rangle\\
&=\pm(\overline{\langle U,S\rangle}-\frac{t^k-1}{m}\overline{\langle U,B\rangle})\\
&=\pm(\overline{V(U,p_*(S))}-\frac{t^k-1}{m}\overline{V(U,p_*(B))})\\
&=\pm(\overline{V(U,(t^k-1)X)}-\frac{t^k-1}{m}\overline{V(U,mX)})\\
&=\pm(t^k-1)(\overline{V(U,X)}-\overline{V(U,X}))=0.
\end{align*}
Note that $(t^k-1)$ gets conjugated twice as it is pulled out of the first term: once due to the Hermitian property of the pairing and once by the explicit conjugation. 

For an alternate choice of $S$, say $S'$, with $p_*(S')=(t^k-1)X$, then $S-S'\in\ker(p_*)$, and choosing $M\in H_n(\td C)$ such that $p_*(M)=mY$, we have $V(S-S',Y)=\frac{1}{m}V(S-S', mY)=\frac{1}{m}V(S-S', p_*(M))=\pm\frac{1}{m}\overline{V(M, p_*(S-S'))}=0$. A similar argument shows independence of choice of $B$. 

Now, if we choose $X'$ instead of $X$ along with corresponding $B'$ and $S'$ and let $M,N\in H_n(\td C)$ such that $p_*(M)=mY$ and $p_*(N)=(t^k-1)Y$, we get 
\begin{align*}
V(S,Y)&-\frac{t^k-1}{m}V(B,Y)-(V(S',Y)-\frac{t^k-1}{m}V(B',Y))\\
 &=V(S-S', Y)-\frac{t^k-1}{m}V(B-B',Y)\\
&=\overline{(t^k-1)}^{-1}V(S-S', (t^k-1)Y)-\frac{t^k-1}{m^2}V(B-B',mY)\\
&=\overline{(t^k-1)}^{-1}V(S-S', p_*(N))-\frac{t^k-1}{m^2}V(B-B',p_*(M))\\
&=\pm(\overline{(t^k-1)}^{-1}\overline{V(N, p_*(S-S'))}-\frac{t^k-1}{m^2}\overline{V(M,p_*(B-B'))}  )\\
&=\pm(\overline{(t^k-1)}^{-1}\overline{V(N, (t^k-1)(X-X'))}-\frac{t^k-1}{m^2}\overline{V(M,m(X-X'))}  )\\
&=\pm(\frac{t^k-1}{\overline{t^k-1}}\overline{V(N, X-X')}- \frac{t^k-1}{m}   \overline{V(M,X-X')}  ).
\end{align*}
Conjugating the whole equation gives
\begin{align*}
\pm(\frac{\overline{t^k-1}}{t^k-1}V(N, X-X')- \frac{\overline{t^k-1}}{m}  V(M,X-X'))\\
\qquad =
\pm\frac{\overline{t^k-1}}{t^k-1}(V(N, X-X')- \frac{t^k-1}{m}  V(M,X-X')).
\end{align*}
But we know that the expression inside the parentheses is $0$ from our calculations above with $Y$ and $Y'$ and noting that this expression is independent of choice of $N$ and $M$ just as in the above proof of independence of $B$ and $S$.

Finally, we show that the isometry class of $(T_{n-1}(\td X), [\,,\,])$ is determined algebraically completely by the isometry class of $(H_n(\td C), \langle \,,\,\rangle)$.  Given an isometry $f:(H_n(\td C),\langle \,,\,\rangle)\to (H_n(\td C'),\langle \,,\,\rangle')$, we obtain a commutative diagram of  exact sequences
\begin{equation*}
\begin{CD}
H_n(\td C)&@>p_*>> & \overline{ \Hom(H_n(\td C), Q(\Lambda)/\Lambda)}& @>>> & H_{n-1}(\td  X)&@>>>&0\\
@V f VV && @V (f^*)^{-1}=h VV && @VV g V\\
H_n(\td C')&@>p'_*>> &  \overline{\Hom(H_n(\td C'), Q(\Lambda)/\Lambda)}& @>>> & H_{n-1}(\td X')&@>>>&0,
\end{CD}
\end{equation*}
in which $p_*$ and $p_*'$ are induced by the $(-1)^{n+1}$-Hermitian pairings $\langle \,,\,\rangle$ and $\langle \,,\,\rangle'$ and  $g$ is induced by $(f^*)^{-1}$ so that the commutativity of the last square is tautological. Here $f^*$ is an isomorphism since $f$ is. To see that the first square commute, we use $V(x,p_*(y))=\langle x,y\rangle=\langle f(x),f(y)\rangle'=V'(f(x), p_*'f(y))=V(x, f^*p_*'f(y))$. Since $x$ and $y$ are arbitrary and the pairings are nondegenerate, we we see that $p_*=f^*p_*'f$. So now by reversing $f^*$ to $(f^*)^{-1}$, we get an isomorphism of  exact sequences which induces the isomorphism $g:H_{n-1}(\td X) \to H_{n-1}(\td  X')$.

So $g$ is an isomorphism of modules and, in particular, the same choices of $m$ and $t^k-1$ serves to annihilate both. Finally, we want to show that $g$ induces an isometry of the Farber-Levine pairing $[\,,\,]$ on $T_{n-1}(\td X)\cong T_{n-1}(\td X')$. Of course it is sufficient to show that it induces an isometry of $\{\,,\,\}$. Given $z,w\in T_{n-1}(\td  X)$, we know $\{g(x),g(y)\}'$ is determined by $V'(S',Y')-\frac{t^k-1}{m}V'(B',Y')$, where $B', S', Y'$ are as defined above. By the commutativity, we can choose $Y$ and $X$ such that $h(Y)=Y'$ and $h(X)=X'$, which then implies that there are $S$ and $B$ such that $f(B)=B'$ and $f(S)=S'$. Let us also choose $M$ such that $p_*(M)=mY$ and a corresponding $M'=f(M)$ (so $p_*(M')=mY'=mh(Y)$). We have 
\begin{align*}
V'(S',Y')-\frac{t^k-1}{m}V'(B',Y')&=V'(f(S),h(Y))-\frac{t^k-1}{m}V'(f(B),h(Y))\\
&=\frac{1}{m}V'(f(S),mh(Y)-\frac{t^k-1}{m^2}V'(f(B),mh(Y))\\
&=\frac{1}{m}V'(f(S),h(p_*(M)))-\frac{t^k-1}{m^2}V'(f(B),hp_*(M))\\
&=\frac{1}{m}V'(f(S),p'_*f(M))-\frac{t^k-1}{m^2}V'(f(B),p'_*f(M))\\
&=\frac{1}{m}\langle f (S),f(M)\rangle'-\frac{t^k-1}{m^2}\langle  f(B),f(M)\rangle'\\
&=\frac{1}{m}\langle  S,M\rangle-\frac{t^k-1}{m^2}\langle B,M\rangle\\
&=\frac{1}{m}V(S,p_*(M))-\frac{t^k-1}{m^2}V(B,p_*(M))\\
&=V(S,Y)-\frac{t^k-1}{m}V(B, Y).
\end{align*}
So $g$ is an isometry of the Farber-Levine pairing.

We summarize what we have achieved so far  as follows:

\begin{theorem}\label{T: Blanch gives FL}
Given a simple disk knot $D^{2n-1}\into D^{2n+1}$, the module $H_{n-1}(\td X)$ and the Farber-Levine $\Z$-torsion pairing on $T_{n-1}(\td X)$  are determined up to isometry by the isometry class of the Blanchfield self-pairing on $H_n(\td C)$. 
\end{theorem}

With a little more work, one could enlarge this theorem to apply to more general cases, for example some disk knots that are not necessarily simple. However, the theorem as stated will be sufficient for our purposes.

Before we move on to showing that $T_{n-1}(\td X)$ and the Farber-Levine pairing really only depend on $\bar H$ and its self-Blanchfield pairing, it is worth pausing to clear up the dependence of our formulation of the Farber-Levine pairing on the choices of $m$ and $k$. Although we began the purely algebraic part of our discussion by fixing a canonical $k$ and $m$, i.e. we chose $k$ and $m$ without any ambiguity as the smallest positive integers satisfying certain properties, let us be complete and demonstrate independence of $k$ and $m$ within the restrictions imposed by these properties. Of course the Farber-Levine pairing does not depend on these choices by the work in \cite{L77}, but we will compute directly from our new definition. 

First, let us consider a new choice $k'$. We need $t^{k'}=1$ on certain modules, and we have assumed that $k$ is the smallest positive integer for which this holds, so we must have $k'=bk$ for some $0<b\in \Z$. Now recall that to define $\{x,y\}$ in terms of $V$, we had to find $X, Y\in H_n(\td C,\td X)$ such that $\bd(X)=x$ and $\bd (Y)=y$ and $B,S\in H_n(\td C)$ such that $p_*(B)=mX$, $p_*(S)=(t^{k}-1)X$. The only choice that depends on $k$ is that of $S$, so let us take a new $S'\in H_n(\td C)$ with $p_*(S')=(t^{k'}-1)X$. In fact, let us choose $S'=\frac{t^{k'}-1}{t^k-1}S$, which will suffice (note that $\frac{t^{k'}-1}{t^k-1}\in \Lambda$). Now we compute
\begin{align*}
V(S',Y)-\frac{t^{k'}-1}{m}V(B,Y) &= V( \frac{t^{k'}-1}{t^k-1}S, Y)-\frac{t^{k'}-1}{t^k-1}\frac{t^k-1}{m}V(B,Y)\\
  &=  \frac{t^{k'}-1}{t^k-1}(V(S,Y)-\frac{t^{k}-1}{m}V(B,Y)).
\end{align*}
If we let $\theta_k=\Lambda/\langle t^k-1 \rangle$, then the final formula gets projected to $I(\theta_k)/\theta_k$ while the first gets projected to $I(\theta_{k'})/\theta_{k'}$. So the question becomes whether the diagram
\begin{diagram}
T_{n-1}(\td X)\otimes T_{n-1}(\td X) &\rTo^{\{\,,\,\}_k} & I(\theta_k)/\theta_k \\
\dTo^{ \{\,,\,\}_{k'}} & \ldTo(2,2)^{\frac{t^{k'}-1}{t^k-1}}& \dTo\\
I(\theta_{k'})/\theta_{k'}  &\rTo &\Q/\Z
\end{diagram}
commutes on the bottom triangle. This follows from purely algebraic computations of Levine \cite[pp. 12, 16]{L77}.

Now let's see what happens if we change $m$ to $m'$. Since $m$ was also selected as the smallest positive integer which kills certain groups, we must have $m'=am$ for some $0<a\in \Z$. Recall that $k'$ must now be chosen so that $t^{k'}=1$ on a certain finite collection of finite $\Lambda$-modules depending on $m'$ and that $k$ is chosen similarly to correspond to $m$. But then there must exist a $K$ such that $t^K=1$ on all of these modules (e.g. just take $K=kk'$). Since we have already shown the pairing to be independent of choice of $k$ for fixed $m$, we are free to work with this $K$. Once again we can choose $X$ and $Y$ independent of $m$ or $k$, and we choose $B$, $B'$, and $S$ so that $p_*(B)=mX$, $p_*(B')=m'X$, and $p_*(S)=(t^K-1)X$. In fact, we can choose $B'=aB$. Then 
\begin{align*}
V(S,Y)-\frac{t^{K}-1}{m'}V(B',Y) &= V(S, Y)-\frac{t^K-1}{am}V(aB,Y)\\
  &= V(S,Y)-\frac{t^{K}-1}{m}V(B,Y).
\end{align*}
And since we have established independence of $K$, this equation gives us independence of $m$ on passage to $\Q/\Z$.

\paragraph{The Farber-Levine pairing depends only on $\mathbf{\bar H}$ and its pairing}

Now that we have shown that, for a simple disk knot, the Farber-Levine torsion pairing on $T_{n-1}(\td X)$ is determined by the self-Blanchfield pairing on $H_n(\td C)$, we wish to strengthen this result somewhat and show that, in fact, it only depends on the self-Blanchfield pairing on $\bar H$, the cokernel of the map $H_n(\td X)\to H_n(\td C)$. This pairing will no longer determine all of $H_{n-1}(\td X)$, but it suffices to determine $T_{n-1}(\td X)$ and its Farber-Levine pairing. From this, we will be able to conclude that the Farber-Levine pairing is determined by the Seifert matrix of the disk knot.

\begin{theorem}\label{H to FL}
For a simple disk knot $L:D^{2n-1}\subset D^{2n+1}$, the $\Lambda$-module $T_{n-1}(\td X)$ and its Farber-Levine $\Z$-torsion pairing are determined up to isometry by the isometry class $\bar H$ with its self-Blanchfield pairing.
\end{theorem}

\begin{proof}

Once again, we know that we have the exact sequence
\begin{equation*}
\begin{CD}
0 & @>>> &  H_n(\td X) & @> i_* >> & H_n(\td C) & @>p_*>> & H_n(\td C, \td X) & @>>> & H_{n-1}(\td X) & @>>> & 0
\end{CD}
\end{equation*}
and that the modules $H_{n-1}(\td X)$ and $T_{n-1}(\td X)$ and the Farber-Levine pairing on $T_{n-1}(\td X)$ are determined by the self-Blanchfield pairing on $H_n(\td C)$.  The module $\bar H$ is the cokernel of $i_*$, and it contains no $\Z$-torsion as $H_n(\td C,\td X)$ is $\Z$-torsion free (since the knot is simple). In the purely algebraic context also considered above, we could alternatively define $\bar H$ as $H_n(\td C)/H^{\perp}_n(\td C)$, where $H^{\perp}_n(\td C)$ is the annihilator of $H_n(\td C)$ under its Blanchfield pairing. In this context, $\bar H$ has no $\Z$-torsion because it injects into $\overline{ \Hom_{\Lambda} (\bar H , Q(\Lambda)/\Lambda)}$, which has no $\Z$-torsion by \cite[Props. 3.2, 4.1]{L77}.  

Consider now the following diagram:

\begin{diagram}
0 &  \rTo & \bar H & \rTo^{p} & \overline{ \Hom_{\Lambda} (H_n(\td C), Q(\Lambda)/\Lambda)} & \rTo^{\bd} & H_{n-1}(\td X) & \rTo & 0 \\
    &           & \uTo^{=}  & & \uInto^{\pi^*}  &  & \uInto^g\\
0 &  \rTo & \bar H & \rTo^{\rho} & \overline{ \Hom_{\Lambda} (\bar H , Q(\Lambda)/\Lambda)} & \rTo^{\eta} & A & \rTo & 0 &.
\end{diagram}
Denoting by $p$ the map induced from $p_*$, the first line is exact by the exactness of the preceding diagram, which also induces the map $\rho: \bar H \to \overline{ \Hom_{\Lambda} (\bar H , Q(\Lambda)/\Lambda)}$ since the self-Blanchfield pairing is trivial on any element of $\im(i_*)=\ker(p_*)=H^{\perp}_n(\td C)$. The map $\rho$ is injective since the self-Blanchfield pairing on $\bar H$ is nondegenerate. 
The map $\pi^*$ is induced by the projection $\pi: H_n(\td C)\to \bar H$, and it is injective since $\pi$ is surjective and the $\Hom$ functor is left exact. The $\Lambda$-module $A$ is the cokernel of $\rho$ by definition, and $g$ is induced by the rest of the diagram. $g$ is injective by the five-lemma. 

Suppose $x\in T_{n-1}(\td X)$, the $\Z$-torsion submodule, and that $m=|T_{n-1}(\td X)|$. From the diagram, $x=\bd (y)$ for some $y\in  \overline{ \Hom_{\Lambda} (H_n(\td C), Q(\Lambda)/\Lambda)}$, and $my=p(z)$ for some $z\in \bar H$. By commutativity, $my=\pi^*\rho(z)$. So $my$ lifts to $ \overline{ \Hom_{\Lambda} (\bar H , Q(\Lambda)/\Lambda)}$, which means that $my$ annihilates the subgroup $H_n(\td X)$ of $H_n(\td C)$. So for every element $w\in H_n(\td X)$, $my(w)=V(w,my)\in \Lambda$, which implies that each rational function $y(w)$ must be of the form $\lambda_{w}/m$ for some $\lambda_{w}\in \Lambda$. 

We claim that in fact we must then  have $y(w)\in \Lambda $ for every $w\in H_n(\td X)$. The proof is similar to that of  \cite[Lemma 5.1]{L77}. First note that in the abstract algebraic context, $H_n(\td X)$ is of type K by applying the five lemma to the above commutative diagram mapped by $t-1$; we  note that $H_n(\td C)$ is of type K by algebraic assumption and  $\overline{ \Hom_{\Lambda} (H_n(\td C), Q(\Lambda)/\Lambda)}$ is of type K by \cite[Prop 4.1 and p. 8]{L77}. Now, suppose that $w\in H_n(\td X)$. Since $H_n(\td X)$ is of type K, by the proof of  \cite[Cor. 1.3]{L77} there is a polynomial $\Delta$ such that $\Delta H_n(\td X)=0$ and $\Delta(1)=\pm1$. So $y(\Delta w)=\Delta y(w)=\Delta \lambda_w/m \in \Lambda$. But since $\Delta(1)=\pm1$, no factor of $m$ divides $\Delta$ in $\Lambda$, so it must be that $m$ divides each $\lambda_w$, i.e. $y(w)=\lambda_w/m\in \Lambda$. 

This shows that $y$ annihilates $H_n(\td X)$, which implies that $y$ lifts to an element in  $\overline{\Hom( \bar H, Q(\Lambda)/\Lambda)}$, i.e. $y=\pi^*(y')$, which implies that $x=\bd(y)=g\eta(y')$. So $T_{n-1}(\td X)\subset \im(g)$. But $T_{n-1}(\td X)$ is finite \cite[Lemma 3.1]{L77} and $g$ is injective, so we must have $T_{n-1}(\td X)\cong T(A)$, the $\Z$-torsion subgroup of $A$. 

Of course by our previous discussion, the Farber-Levine $\Z$-torsion pairing on $T(A)$ is determined by the self-Blanchfield pairing on $\bar H$, and by the inclusion of the second row of the diagram into the first, these parings are compatible with those of $T_{n-1}(\td X)$ and $\bar H$ (as induced from $H_n(\td C)$). So $T(A)$ and $T_{n-1}(\td X)$ are isomorphic with isometric Farber-Levine pairings, induced by the self-Blanchfield pairing on $\bar H$. 

\end{proof}

\begin{corollary}\label{C: Seifert gives FL}
For a simple disk knot $L:D^{2n-1}\subset D^{2n+1}$, the $\Lambda$-module $T_{n-1}(\td X)$ and its Farber-Levine $\Z$-torsion pairing are determined up to isometry by any Seifert matrix for $L$.
\end{corollary}
\begin{proof}
As seen in Section \ref{S: realize}, the module $\bar H$ and its Blanchfield self-pairing are determined by any Seifert matrix for $L$. Thus the corollary follows immediately from the preceding theorem. 
\end{proof}

\section{Changing Seifert surfaces}\label{S: seifert surfaces}

The entirety of this long section will be devoted to studying what happens to the Seifert matrix of a disk knot when we change the Seifert surface. Such alterations can always be performed by first doing surgery on the boundary Seifert surface $F$ and then performing internal surgeries that avoid the boundary. Although we will see that different effects arise in different cases, we can summarize the results as follows:

\begin{theorem}\label{T: s-equiv}
Any two Seifert matrices for a disk knot differ by a rational S-equivalence. 
\end{theorem}

Proposition \ref{L: s-equivalence cobordism}, which stated that two Seifert matrices for a disk knot are  cobordant, follows. 

To prove the theorem, we need to relate various Seifert surfaces for a fixed disk knot. So suppose that we have two copies of a disk knot $L$, which we will call $L_1$ and $L_2$, with Seifert surfaces $V_1$ and $V_2$ and boundary Seifert surfaces $F_1$ and $F_2$. Consider the knot $L\times I =(D^{2n+1}\times I, D^{2n-1}\times I)$. This is also a disk knot, and we can think of it as realizing the trivial cobordism from $L_1$ to $L_2$. On the boundary, $\bd D^{2n+1}\times I$, we have the trivial cobordism of the boundary knot $K$. As in \cite[\S]{L70}, we can then construct a cobordism $U$ from $F_1$ to $F_2$ in $\bd D^N\times I$ such that $\bd U$ is equal to the union of $F_1$, $-F_2$ and the trace of the trivial isotopy. The union $V_1\cup U\cup -V_2$ is a Seifert surface for $L_1\cup K\times I\cup  -L_2$, the boundary knot of $L\times I$. By \cite[\S 8]{L69}, this can we extended to a Seifert surface $W$ for $L\times I$. The pair $(W,U)$ thus provides a cobordism from $(V_1,F_1)$ to $(V_2,F_2)$. 

Now, as usual when dealing with cobordism with boundaries, we can break up the process into two distinct steps. We can first consider the cobordism of the boundary. In our case this amounts to beginning with $V_1\subset D^{2n+1}$ and adjoining $U\subset \bd D^{2n+1}\times I$. In other words, we form $(D^{2n+1},V_1)\cup_{(S^{2n}, F_1)}(S^{2n}\times I, U)$. 
Note that we do not need to mention the knots explicitly since they are contained in the embedding information. Then we perform the usual trick and ``rekink'' the diagram so that $W$ becomes a cobordism rel boundary from $V_1\cup U$ to $V_2$. 

In the first subsection below, we consider the second stage and determine how a Seifert matrix is affected by an internal  cobordism, i.e. one that leave the boundary Seifert surface fixed. In the second subsection, we consider the effect of the boundary cobordism.

\subsection{Changing the Seifert surface on the interior}

In this subsection, we first assume that we have two of the same disk knot $L: D^{2n-1}\into D^{2n+1}$ (denoted $L_1$ and $L_2$ when necessary) with two Seifert surfaces $V_1$ and $V_2$ that agree on the boundary (i.e. they have the same Seifert surface for the boundary sphere knot $F:=F_1=F_2$), then we can embed $L\times I$ in $D^{2n+1}\times I$ and consider the boundary knot $L_1\cup -L_2\cup K\times I$ and its Seifert surface $V_1\cup -V_2\cup F\times I$. This can be extended to a Seifert surface $W$ for the whole disk knot $L\times I$, see \cite[\S 8]{L69}. We can now proceed to analyze the change in the Seifert matrix from that obtained from $V_1$ to that obtained from $V_2$ analogously to the case for sphere knots in \cite{L70}. In particular, we can create a smooth (PL) height function and separate $W$ into critical levels. This allows us to restrict to the case where $W$ is obtained from $V_1$ by adding a single handle so that $V_1$ and $V_2$ differ by a single surgery. We will make this assumption throughout.

As in \cite[Lemma 1]{L70}, if cycles in $V_1$ and $V_2$ are homologous through $W$, then they admit the same linking pairing numbers. We state this as a lemma for future reference: 

\begin{lemma}\label{L: Levine lemma}
If $\alpha,\beta\in F_n(V_1)$, $\alpha',\beta'\in F_n(V_2)$, and $i_1(\alpha)=i_2(\alpha')$, $i_1(\beta)=i_2(\beta')$, where $i_j: F_n(V_j)\to F_n(W)$ are induced by inclusion, then $\theta_1(\alpha,\beta)=\theta_2(\alpha',\beta')$, where $\theta$ and $\theta_2$ are the Seifert pairings on $V_1$ and $V_2$, respectively. The analogous results holds for the Seifert pairings relating $F_n(V_j)$ and $F_n(V_j,F)$.
\end{lemma}

\begin{proof}
This is an immediate generalization of \cite[Lemma 1]{L70}.
\end{proof}

Also as in \cite{L70}, there is no effect to  $F_n(V_1)$  if the index of the handle is less than $n$ or if it has index $n$ and the boundary of the cocore of the handle (which is $\bd _*$ of a generator of $H_{n+1}(W, V_2)$) has finite order in $H_n(V_2)$: If the index is $i<n$, then $H_*(W,V_1)$ is $\Z$ for $*=i$ and $0$ otherwise, while $H_*(W,V_2)$ is $\Z$ for $*=2n-i$ and $0$ otherwise. It follows from the long exact sequences of the pairs that $H_n(V_1)\cong H_n(W)\cong H_n(V_2)$. We also have a commutative diagram induced by inclusions

\begin{equation}\label{E: homologies}
\begin{CD}
H_n(F\times 1) & @>>> & H_n(V_2)\\
@V\cong VV & & @VV\cong V\\
H_n(F\times I) & @>>> & H_n(W)\\
@A\cong AA & & @AA\cong A\\
H_n(F\times 0) & @>>> & H_n(V_1) &\,\,,
\end{CD}
\end{equation}
so $\bar E_1\cong \bar E_2$ and this is an isometry of Seifert pairings by Lemma \ref{L: Levine lemma}. (Recall our notation from Section \ref{S: basics}.)

If the index of the handle is $n$ but the boundary of the cocore represents a torsion element $a\in H_n(V_2)$, we consider the following diagram

\begin{diagram}[LaTeXeqno]\label{D: Levine's diagram}
&                    &                   &         &   H_{n+1}(W, V_2)\cong \Z\\
  &                  &                   &         &    \dTo\\
    &                &                   &         &   H_n(V_2)\\
      &              &                   &         &  \dTo\\
0&\rTo & H_n(V_1) &\rTo& H_n(W) &\rTo^{\bd_*} &H_n(W, V_1)\cong \Z\\
    &                &                   &         & \dTo\\
&&&&0&.
\end{diagram}

\noindent The composition $H_n(V_2)\to H_n(W)\to H_n(W, V_1)$ is determined by the intersection number in $V_2$ of the cycle in $H_n(V_2)$ with that given by the boundary of the cocore of the handle (see \cite{L70}; roughly, this intersection number measures how many times a cycle wraps around the handle). If this boundary $a$ has finite order, then this composition is $0$ since it maps to $\Z$. Hence $H_n(V_2)$ maps into the the kernel of $\bd_*$, which is the injective image of $H_n(V_1)$. But  this map is also onto $H_n(W)$ with kernel a torsion element. So $H_n(W)\cong H_n(V_1)$, and these are isomorphic to $H_n(V_2)$ modulo the torsion subgroup generated by $a$. 
Since diagram \eqref{E: homologies} continues to hold modulo torsion, $\bar E_1\cong \bar E_2$, and we can pick  a basis for $F_n(V_2)$ corresponding to our fixed one for $F_n(V_1)$ through homologies in $W$. So by Lemma \ref{L: Levine lemma}, $\bar E_1$ and $ \bar E_2$ are isometric with respect to the Seifert pairing.

\vskip1cm

Suppose now that we have a handle of index $n$ and that the boundary of the cocore, $a\in H_n(V_2)$, is not a torsion element and that $a_0$ is its primitive (i.e. $a$ is a non-trivial  positive multiple of $a_0$ and $a_0$ is not a multiple of any other element).  We claim that either $\bar E$ and its pairing are unaffected or that  $F_n(V_2)\cong F_n(V_1)\oplus \Z \oplus \Z$ and  $\bar E\cong \bar E\oplus \Z^2$. Note that this claim says nothing yet about the extension of the pairing; we shall discuss this below.

Consider again the above commutative diagram \eqref{D: Levine's diagram}, but with all homology groups replaced by those with rational coefficients. By assumption, the map $\Q\cong H_{n+1}(W,V_2;\Q)\to H_n(V_2;\Q)$ is injective, which implies that $H_n(V_2;\Q)\cong H_n(W;\Q)\oplus \Q$. Suppose that $H_n(W,V_1;\Q)\to H_{n-1}(V_1; \Q)$ is non-trivial and hence injective. Then $H_n(V_1;\Q)\cong H_n(W;\Q)$ and $H_n(V_2;Q)\cong H_n(V_1; \Q)\oplus \Q$. We can leave the basis of $F_n(V_1)$ as the basis of $H_n(V_1;\Q)$ and translate via homologies in $W$ to basis elements of $H_n(V_2;\Q)$. Then a basis of $H_n(V_2;\Q)$ can consist of these together with $a$. But since $a$ is the boundary of the cocore of the handle, it is clear that, in $V_2$, $a$ intersects trivially both  itself and all the generators from $H_n(V_1)$, which can be taken to lie in $V_0=V_1-S^{n-1}\times D^{n+1}$ ($V_1$ minus the attaching sphere) by general position. Thus the matrix $R_2$ representing $H_n(V_2;\Q)\to H_n(V_2,F;\Q)$ would increase by a row and column of $0$s from $R_1$; this shows that $a$ must generate a summand in the kernel of $p_*$. Thus $\bar E$ is unaffected, and the linking pairings of the remaining basis elements remain unchanged by Lemma \ref{L: Levine lemma}. So  $\bar E$ remains unchanged up to  isometry. 

So, now suppose that $H_n(W,V_1;\Q)\to H_{n-1}(V_1;\Q)$ is not injective, which means that it must  be $0$ rationally. Then we see that $H_n(W;\Q)\cong H_n(V_1;\Q)\oplus \Q$, which implies that $H_n(V_2;\Q)\cong H_n(V_1)\oplus \Q^2$. This implies that $F_n(V_2)\cong F_n(V_1)\oplus \Z^2$, and we further claim that $\bar E_2\cong \bar E_1\oplus \Z^2$. This will be accomplished if we show that $\bar E_2 \otimes \Q \cong (\bar E_1\otimes \Q )\oplus \Q^2$.  To see this, we will make some more calculations.

Let $V_0=V_1-S^{n-1}\times D^{n+1}$ be $V_1$ minus a neighborhood of attaching sphere.  $H_i(V_0)\cong H_i(V_1)$ for $i\leq n-1$ by general position. 
Next, we consider the Mayer-Vietoris sequence 

\begin{equation*}
\begin{CD}
@>>>& H_i(S^{n-1}\times S^n)&@>k_*>>& H_i(V_0)\oplus H_i(D^{n}\times S^{n}) &@>\rho_* >>& H_i(V_2)&@>>>.
\end{CD}
\end{equation*}
We see that $H_i(V_0)\to H_i(V_2)$ is an isomorphism induced by inclusion for $i\leq n-2$. For $i=n-1$, it is  a rational isomorphism: The map $k_*$ in dimension $n$ from $H_n(S^{n-1}\times S^n)\cong \Z$ to $H_n(V_0)\oplus H_n(D^n\times S^n)$ is injective since the generator of $H_n(S^{n-1}\times S^n)$ is also a generator of the summand $H_n(D^n\times S^n)$. So $H_n(V_2)$ is isomorphic to the direct sum of the image of $\rho_*$, which has rank equal to that of $H_n(V_0)$, and at most one $\Z$ summand which would come from  $H_{n-1}(S^{n-1}\times S^n)$. But from the Mayer-Vietoris sequence for $V_1$,
\begin{equation*}
\begin{CD}
&@>>>& H_i(S^{n-1}\times S^n)&@>>>& H_i(V_0)\oplus H_i(S^{n-1}\times D^{n+1}) &@>j_*>>& H_i(V_1) &@>>>&,
\end{CD}\end{equation*} 
we know that $H_n(V_0;\Q)$ has at most $1$ dimension more than $H_n(V_1;\Q)$. This is because $H_{n-1}(S^{n-1}\times S^n)\to  H_{n-1}(V_0)\oplus H_{n-1}(S^{n-1}\times D^{n+1})$ must be injective as it carries a generator of $H_{n-1}(S^{n-1}\times S^n)$ to a generator of $H_{n-1}(S^{n-1}\times D^{n+1})$. 
So to get the extra two rational summands that we must have in $H_n(V_2;\Q)$ from $H_n(V_1;\Q)$, it must in fact be the case that $k_*\otimes \Q$ is $0$ in dimension $n-1$ (and also then that $F_n(V_0)\cong F_n(V_1)\oplus \Z$ with the projection $F_n(V_0)\to F_n(V_1)$ induced by the inclusion map).

So now consider the diagram
\begin{diagram}[LaTeXeqno]\label{D: surgeries}
H_n(V_2)&\rTo^{r_2} &H_n(V_2, F)&\rTo^{\bd _2} &H_{n-1}(F)&\rTo^{i_2} & H_{n-1}(V_2)&\rTo & H_{n-1}(V_2,F)\\
\uTo^{\theta_2}& &\uTo^{\phi_2}& &\uTo^{=}& &\uTo^{\eta_2}& &\uTo& &\\
H_n(V_0)&\rTo^{r_0} &H_n(V_0, F)&\rTo &H_{n-1}(F)&\rTo & H_{n-1}(V_0)&\rTo & H_{n-1}(V_0,F)\\
\dTo^{\theta_1}& &\dTo^{\phi_1}& &\dTo^{=}& &\dTo^{\eta_1}& &\dTo& &\\
H_n(V_1)&\rTo^{r_1} &H_n(V_1, F)&\rTo^{\bd_1} &H_{n-1}(F)&\rTo^{i_1} & H_{n-1}(V_1)&\rTo & H_{n-1}(V_1,F)&.
\end{diagram}
We know that all vertical maps from $H_{n-1}(F)$ and to the right are \emph{rational} isomorphisms (using the five lemma for the relative terms). We also know by Poincar\'e-Lefschetz duality that $H_n(V_2,F;\Q)$ has dimension two greater than $H_n(V_1, F;\Q)$, since this is true of $H_n(V_2;\Q)$ and $H_n(V_1;\Q)$. But the $\Q$-dimension of $\ker (i_1)$ is the same as that of $\ker (i_2)$ according to the diagram. Thus the $\Q$-dimensions of $\im (\bd _2)$ and $\im (\bd_1)$ agree. But these are equal to the dimensions of the coimages of these boundary maps, which implies that the kernel of $\bd_2$ must have dimension two greater than that of $\bd_1$. Hence we see that $E_2\otimes \Q\cong (E_1\otimes \Q)\oplus \Q^2$. Hence $\bar E_2\otimes \Q\cong (\bar E_1\otimes \Q)\oplus \Q^2$ since $E$ and $\bar E$ always have the same rank. Thus the desired correspondence holds for the ranks of the integral groups.

If we look at diagram \eqref{D: surgeries} integrally, the map $\eta_1$ remains an isomorphism, as observed above, but the map $\eta_2$ may only be onto with kernel a torsion subgroup. It is onto by the Mayer-Vietoris sequence. That the kernel is torsion follows since we have seen that the kernel is rationally $0$. We also know that $\theta_1$ and $\phi_1$ are onto by general position.  

Now, let's think about $\theta_2$. We've established that in the Mayer-Vietoris sequence for $V_2$ that the image of $H_n(V_2)\to H_{n-1}(S^{n-1}\times S^n)$ must be isomorphic to $\Z$. Therefore, we have an exact sequence 
\begin{equation*}\label{E: V2 split}
\begin{CD}
0&@>>>&   [H_n(V_0)\oplus H_n(D^n\times S^n)]/\im(k_*) &@>>>&H_n(V_2) &@>>> &\Z &@>>>0.
\end{CD}
\end{equation*}
But again the image of $k_*$ is of the form $(z,1)$ since a generator of $H_n(S^{n-1}\times S^n)$ maps to a generator of $H_n(D^n\times S^n)$, and furthermore, this sequence must split since it is onto a $\Z$ term. Thus  the map $\theta_2$ is split injective, and, modulo torsion, $H_n(V_2)\cong H_n(V_0)\oplus \Z\cong H_n(V_1)\oplus \Z^2$. One of the $\Z$ terms is generated by the primitive of the cocore of the handle, which was $a_0$; in fact $F_n(V_0)\cong F_n(V_1)\oplus \langle a_0\rangle$ from the Mayer-Vietoris sequence.

As for $\phi_2$, it must also then be injective by the five lemma.

So at this point, we have established split exact sequences 
\begin{equation*}
\begin{CD}
0& @>>> & \Z & @>>> & F_n(V_0) & @>>> & F_n(V_1) & @>>> & 0\\
0& @>>> & H_n(V_0) & @>>> & H_n(V_2) & @>>> & \Z & @>>> & 0&\,\,.
\end{CD}
\end{equation*}
Now recall from Section \ref{S: basics} that we can choose a basis $\{ \delta_i\}_{i=1}^k$ for $F_n(V_1)$ such that $\{ \delta_i\}_{i=1}^m$ project to a basis $\{\bar  \delta_i\}_{i=1}^m$ for $\bar E_1$ and $\{ \delta_i\}_{i=m+1}^k$ are a basis for the kernel of $r_1: F_n(V_1)\to F_n(V_1,F)$. Since the kernel of $r_1$ is equal to the image of the map induced by the inclusion $F\to V_1$, we can further assume that these latter classes are represented by chains lying in $F$. By general position, we can also assume that we have chosen chains representing the $\delta_i$ that do not intersect a neighborhood of the attaching sphere of the surgery. This chooses a lift of these basis elements to $F_n(V_0)$ and using the first split exact sequence above, we see that a basis of $F_n(V_0)$ consists of these lifts plus an element generating the extra $\Z$ summand. And from the Mayer-Vietoris sequence for $V_1$, the extra $\Z$ summand is generated by a primitive of the element represented by $*\times S^n$, which in $V_2$ will be homologous to the boundary of the cocore. In fact, because of the splitting, this primitive will be the boundary of the cocore itself.  So now using the split injection of the second short exact sequence, we see that a basis for $F_n(V_2)$ consists of the lifted $\delta_i$, the boundary of the cocore which we have called $a_0$, and a generator of the second new $\Z$ summand, which, we will call $b_0$ and which maps nontrivially into $H_{n-1}(S^{n-1}\times S^n)$ in the Mayer-Vietoris sequence for $V_2$. Furthermore, since the $\{ \delta_i\}_{i=m+1}^k$ are still represented by chains in $F$, they are also in the kernel of $r_2$. But by our earlier dimension counting arguments, which tell us that we must have $\bar E_2\cong \bar E_1\oplus \Z^2$, this implies that this is the entire kernel of $r_2$, so $\bar E_2$ has a basis consisting of the $\{ \delta_i\}_{i=1}^m$, $a_0$, and $b_0$.

So, at the moment, we have a basis of $\bar E$ consisting of the images of the $\{\delta_i\}_{i=1}^m$ (translated up to $V_2$) plus  $a_0$ and $b_0$. 
Now we can finally look at the Seifert matrix for $\bar E_2$. Lemma \ref{L: Levine lemma} still holds in that these translated  $ \delta_i$ will have the same linking matrix $\theta_1$ as they did for $V_1$. And since $a_0$ corresponds to the boundary of the cocore, it is null-homologous in $W$ and thus links trivially with all the translated $\delta_i$ and also with itself.
Thus as in \cite{L70}, we obtain a matrix for $\theta_2$ of the following form:
\begin{equation}
\begin{pmatrix}\label{E: S-equiv}
\theta_1 & 0 & \eta\\
0 & 0 &x\\
\xi & x' & y
\end{pmatrix},
\end{equation}
where $\theta_1$ is an $m\times m$ matrix, $\eta$ is a $1\times m$ matrix, $\xi$ is an $m\times 1$ matrix, and $x$, $x'$, and $y$ are integers. The key difference from Levine's matrix \cite[p. 188]{L70} is that here $x'+(-1)^nx$, while being the intersection number of $a_0$ and $b_0$, will not necessarily be $\pm 1$. 

In fact, this element $b_0$ maps to a multiple of the generator of $H_{n-1}(S^{n-1}\times S^n)$ under the boundary map of the Mayer-Vietoris sequence for $V_2$. This implies that as a chain, $b_0$ can be represented by a multiple of the attached disk $D^n$ suitably translated into $V_2$ plus another piece whose boundary is a multiple of the attaching $S^{n-1}$, also translated into $V_2$. Note that the intersection number of $b_0$ and $a_0$ is the smallest possible (in absolute value) non-zero intersection number between $a_0$ and all elements of $\bar E_2$:  $a_0$ does not intersect any of the  $\delta_i$, since they all lie in $V_0$ and  $a_0$ is the cocore of the handle. Nor does $a_0$ intersect itself, since the cocore can be pushed off itself along the handle. So no further changes of basis keeping $a_0$ fixed can provide a basis element that has a smaller non-zero intersection number with $a_0$ than $b_0$ does. Clearly, however, the intersection of $a_0$ and $b_0$ is non-trivial.

Now, from \cite[\S 3.6]{GBF1}, the Alexander polynomial $c_n(t)$ associated to the coimage of $H_n(\td C;\Q)\to H_n(\td C, \td X;\Q)$  and determined up to similarity in $\Lambda$ is the determinant of $(-1)^{n+1}(R^{-1})'\tau Rt-\tau'=(R^{-1})'((-1)^{n+1}\tau R t- R'\tau')=(R^{-1})'((-1)^{n+1}\theta' t- \theta)$. But recall that we also know that, with an appropriate integrally unimodular  change of bases (which therefore won't affect its determinant), $-R=\theta+(-1)^n\theta'$, where here $R$ is just the transpose of the intersection matrix on $\bar E$. So the Alexander polynomial is the product of the determinants  of  $((-1)^{n+1}\theta'-\theta))^{-1}$ and $(-1)^{n+1}\theta' t- \theta$. 
If we compare these polynomials as obtained using $\theta_2$ and $\theta_1$, we see that, just as in \cite{L70},  the determinant of $((-1)^{n+1}\theta_2' t- \theta_2)$ is that of $((-1)^{n+1}\theta_1' t- \theta_1)$ multiplied by $((-1)^{n+1}xt-x')((-1)^{n+1}x't-x)$, and we also see that the determinant of $((-1)^{n+1}\theta_2'-\theta_2))^{-1}$  is that of $((-1)^{n+1}\theta_1'-\theta_1)^{-1}$ multiplied by $ ((-1)^{n+1}x-x')((-1)^{n+1}x'-x)$. Since this modification to the Seifert matrix cannot change the polynomial, which is an invariant of the knot, beyond multiplication by $\pm$ a power of $t$, it follows that either $x'$ or $x$ must be $0$. 

If it so happens that  $x'+(-1)^nx=\pm 1$, then $\theta_2$ and $\theta_1$ are integrally S-equivalent as in \cite{L70}. In some cases, this will be guaranteed. For example, if the attaching sphere $S^{n-1}$ is nullhomologous in $V_2$, then $b_0$ can be chosen so that the intersection of $a_0$ and $b_0$  is equal to $1$. We already know that $S^{n-1}$ cannot represent a free element of $V_1$, or else  $\bd_*:\Q\cong H_n(W, V_1;\Q)\to H_{n-1}(V_1; \Q)$ will be injective, which will imply that $H_n(W;\Q)\cong H_n(V_1;\Q)$, which we know does not happen in the case under consideration. So the remaining case is that in which $S^{n-1}$ is division null-homologous, but not null-homologous itself. 

We know by Poincar\'e duality that there must be an element of $H_n(V_2, F)$ whose intersection with $a_0$ must be $1$, and again this must be an element that is the sum of two chains, one of which is represented by the core of the handle (pushed into the boundary of the handle) and the other of which must have as boundary one piece that is the attaching sphere and another piece that is in $F$ (this second piece cannot be empty, else $S^{n-1}$ bounds in $V_2$, which is not true in the case under consideration).  In other words, we see that in this case the attaching sphere must be homologous to a cycle in $F$. Thus this ``bad'' case, in which $x'+(-1)^nx\neq\pm 1$,  can only happen if the attaching sphere represents a torsion element of $H_{n-1}(V_1)$ that is in the image of $H_{n-1}(F)$ under inclusion.  In this case, we do not have S-equivalence, per se, but we do obtain a special type of elementary expansion of the form above, with either $x$ or $x'$ equal to $0$ and the other equal to the intersection number of $a_0$ and $b_0$. We do obtain rational S-equivalence

This completes our study of what happens to the Seifert matrix when a handle of index $\leq n$ is added to the interior of $V$. But of course the addition of handles of higher index can be treated by reversing the direction of the cobordism. So this takes care of all surgeries on  spheres in the interior of $V$.

\subsection{Changing the  boundary Seifert surface}

 We have already examined internal surgeries, so it remains to  consider those that simply add to the boundary. Again we can break the situation into the addition of one handle at a time by the usual Morse theory argument. So we must see the effect on the Seifert matrix of adding a handle to $V$ along $F$. We will denote $V$ plus this handle as $V'$, we will let $F'$ be the new resulting boundary piece after the surgery, and we will let $F_0$ represent $F$ minus a neighborhood of the attaching sphere. 

We first prove that in most dimensions attaching a disk to $V$ along $F$ does not affect the Seifert matrix. 

\subsubsection[Handles of index $\neq n,n+1$]{Handles of index $\mathbf{\neq n,n+1}$}
We consider attaching a handle of index $j$ so that $V'\sim_{h.e.}V\cup D^j$. Then $H_i(V',V)\neq 0$ if and only if $i=j$, so $H_i(V)\cong H_i(V')$ by inclusion for $i\neq j, j-1$. In particular, $H_n(V)\cong H_n(V')$ unless $j=n$ or $j=n+1$. 

Meanwhile, we have the Mayer-Vietoris sequences

\begin{equation*}
\begin{CD}
@>>>&H_n(S^{j-1}\times S^{2n-j-1})&@>>>&H_n(F_0)\oplus H_n(S^{j-1}\times D^{2n-j})&@>>>& H_n(F) @>>>\\
@>>>&H_n(S^{j-1}\times S^{2n-j-1})&@>>>&H_n(F_0)\oplus H_n(D^{j}\times S^{2n-j-1})&@>>>& H_n(F') @>>>.
\end{CD}
\end{equation*}
If $j<n-1$ or if $j>n+1$, then we see from these sequences that $H_n(F_0)$ is isomorphic to $H_n(F)$ and $H_n(F')$, the isomorphisms induced by inclusions. So from the diagram

\begin{equation*}
\begin{CD}
H_n(F)&@>i>>&H_n(V)&@>>>& H_n(V,F) @>>>\\
@A k AA && @AAA && @AAA\\
H_n(F_0)&@>i_0>>&H_n(V)&@>>>& H_n(V,F_0) @>>>\\
@VVV && @VVV && @VVV\\
H_n(F')&@>i'>>&H_n(V')&@>>>& H_n(V',F') @>>>,
\end{CD}
\end{equation*}
we also see that $\cok (i)\cong \cok (i')$, and since the maps are induced by inclusions, we can choose the same chains to represent bases of each. The Seifert matrices therefore remain identical, since we then see that we can pick these representative chains in $V$, and their push-offs along normal vector fields to $V$ and subsequent linking numbers are unaltered by the handle addition. 

These arguments can be extended without great difficulty to the case $j=n-1$. The map induced by inclusion $H_n(F_0)\to H_n(F)$ may now fail to be an isomorphism (see the Mayer-Vietoris sequence), but it remains onto, from which it follows that $\cok (i)\cong \cok(i_0)$ since $\im(i_0)=\im(ik)=\im(i)$. The map $H_n(F_0)\to \coim(i')$ must also be onto since we see from the Mayer-Vietoris sequence that there is a surjection $H_n(F_0)\oplus H_n(D^{j}\times S^{2n-j-1})\to H_n(F')$. But the generator of $H_n(D^{j}\times S^{2n-j-1})\cong \Z$ is the boundary of the cocore of the handle and so bounds in $V'$. Thus the image of this summand is in the kernel of $i'$. So $H_n(F_0)$ must map onto the coimage of $i'$ under $H_n(F_0)\to H_n(F')$ followed by projection. Hence the image of $i'$ is equal to the image of the composition of $i'$ with $H_n(F_0)\to H_n(F')$. It now follows again (since $H_n(V)\cong H_n(V')$) that $\cok (i_0)=\cok(i')$ as above but factoring through coimages and the Seifert matrix again remains unchanged as we can choose representative chains in $V$.

This leaves the cases of $j=n$ and $j=n+1$.

\subsubsection[Handles of index $n$]{Handles of index $\mathbf{n}$}

In this case $H_n(V',V)\cong \Z$ and $H_i(V',V)=0$ otherwise. This implies that $H_n(V)\to H_n(V')$ is injective, and either it is an isomorphism or the inclusion of a direct summand, the other summand being $\Z$. 

\paragraph{Case: $\mathbf{H_n(V)\protect\cong H_n(V')}$.} Assume that $H_n(V)\cong H_n(V')$. This will be the case if $\bd_*:H_n(V',V)\to H_{n-1}(V)$ is injective, which will happen if the attaching sphere for the handle generates a free subgroup of $H_{n-1}(V)$. 

The Mayer-Vietoris sequences for $F$ and $F'$ become 

\begin{equation}\label{E: MVF}
\begin{CD}
0&@>>>& H_n(F_0) &@>>>& H_n(F) &@>>>& \Z\oplus \Z &@>\Phi>>& H_{n-1}(F_0)\oplus \Z &@>>>& H_{n-1}(F)&@>>>& 0\\
0&@>>>& H_n(F_0) &@>>>& H_n(F') &@>>>& \Z\oplus \Z &@>>>& H_{n-1}(F_0)\oplus \Z &@>>>& H_{n-1}(F')&@>>>& 0&.\\
\end{CD}
\end{equation}

Claim: the vertical maps induced by inclusion in the following commutative diagram are isomorphisms
\begin{equation*}
\begin{CD}
H_n(F_0)&@>i_0>>& H_n(V)\\
@VVV &&@VVV\\
H_n(F') &@>i'>>& H_n(V'),
\end{CD}
\end{equation*}
and therefore $\cok(i_0)\cong \cok (i')$. The righthand map is an isomorphism by the assumption of this case. The lefthand map is injective by the Mayer-Vietoris sequence. To see that the lefthand map is surjective, we consider the long exact sequence of $(F',F_0)$. By excision, $H_n(F',F_0)\cong H_n(D^n\times S^{n-1}, S^{n-1}\times S^{n-1})\cong \Z$, generated by the cell $D^n\times *$, which is a translate of the attached disk. Then in the exact sequence 
\begin{equation*}
\begin{CD}
H_n(F_0)&@>>>& H_n(F')&@>>>& H_n(F',F_0)\cong \Z&@>\bd_*>>& H_{n-1}(F_0),
\end{CD}
\end{equation*} 
the image under $\bd_*$ of the generator of $H_n(F',F_0)$ is a translate of the attaching sphere in $F_0$. But by the assumptions of this case, we know that this chain generates an infinite cyclic subgroup under the inclusion map $H_{n-1}(F_0)\to H_{n-1}(V)$. Hence $\bd_*$ must be injective, whence $H_n(F_0)\to  H_n(F')$ is surjective. 

We next consider the exact sequence of the pair $(F,F_0)$. By excision, $H_i(F,F_0)\cong H_i(S^{n-1}\times D^n, S^{n-1}\times S^{n-1})$. So again $H_n(F_0)\to H_n(F)$ is injective, and $H_n(F,F_0)\cong \Z$ is generated by the cell $*\times D^n$.

\subparagraph{Subcase: $\mathbf{H_n(F_0)\twoheadrightarrow H_n(F)}$.}  $H_n(F_0)\to H_n(F)$ will be surjective if the boundary of the cell $*\times D^n$, the boundary of a fiber of the normal disk bundle of the attaching sphere, generates an infinite cyclic group in $H_{n-1}(F_0)$. In this case, both vertical maps in 
\begin{equation*}
\begin{CD}
H_n(F)&@>i>>& H_n(V)\\
@AAA &&@AAA\\
H_n(F_0) &@>i_0>>& H_n(V)
\end{CD}
\end{equation*}
are isomorphisms, so $\cok(i_0)\cong \cok (i)$. Thus together with the previous calculation that $\cok(i_0)\cong \cok (i')$, we have $\cok(i)\cong \cok(i')$, and since all of these vertical maps have been by inclusions, each cokernel  can employ the same chains as generators, whence the Seifert matrices are identical. 

\subparagraph{Subcase: Not $\mathbf{H_n(F_0)\twoheadrightarrow H_n(F)}$.} In the alternative case in which a multiple of this fiber sphere bounds in $F_0$, there is a splitting and $H_n(F)\cong H_n(F_0)\oplus \Z$. The $\Z$ term can be generated by the sum of two chains, one lying in $F_0$ and one in $S^{n-1}\times D^n$, both of whose boundary chains are corresponding (opposite sign) multiples of the fiber sphere (of course the one not in $F_0$ will just be a multiple of the fiber disk). This can also be seen from the Mayer-Vietoris sequence. Call this generator $a$. If $a$ bounds in $V$, then $H_n(F_0)\to H_n(F)$ will be onto the coimage of $H_n(F)\to H_n(V)$ and it will follow again that $\cok(i)\cong\cok(i')$. Similarly, if the image of $a$ in $H_n(V)$ is torsion, then $H_n(F_0)\to H_n(F)$ will be onto the coimage of $H_n(F)\to H_n(V)$ mod torsion. Again we get $\cok(i)\cong \cok(i')$ and isometric pairings. 

So the one  remaining case of interest in this subcase will be that in which the image of $a$ generates an infinite cyclic group in $V$. Note that, since $H_n(F_0)\to H_n(F)$ is injective, $H_n(F_0)\to H_n(V)$ actually factors through $H_n(F)$ so that the image of the $H_n(F_0)$ summand of $H_n(F)$ will agree with the image of $H_n(F_0)$. 

We will actually see that a multiple of  the image of $a$ in $H_n(V)$ lies in the image of $H_n(F_0)$. This will imply that  $\cok(i)\cong\cok(i_0)$ mod torsion, and it will follow that the Seifert matrix is unchanged by the addition of the handle. To prove the claim, we consider the image of $a$ in $H_n(V)$, still represented by the chain $a$ as described above. Since the inclusion $H_n(V)\cong H_n(V')$ is an isomorphism, $a$ must represent an infinite cyclic subgroup of $H_n(V')$. The image of this homology class in $H_n(V', F')$, also represented by (the appropriate coset of) $a$, must be $0$ for the following reason. By duality, we know that $H_n(V')$ and $H_n(V',F')$ are dually paired by the intersection form. But our chain representing $a$ in $H_n(V', F')$ can be made disjoint from any other chain representing a class in $H_n(V')$ since all such classes can be assumed to lie in $V$ and hence the interior of $V$ using the inclusion-induced isomorphism $H_n(V)\cong H_n(V')$ and by pushing in along a collar of the boundary $F$ of $V$. But $a$ lies in $F$ and hence is disjoint from any such chain. We conclude that $a$ represents a torsion element in $H_n(V',F')$. Thus some multiple of  $a$ must be in the image of $H_n(F')\to H_n(V')$, and hence the image of the composite $H_n(F_0)\overset{\cong}{\to} H_n(F')\to H_n(V')$. So some multiple of $a$ is representable by a chain lying entirely in $F_0$. By these geometric arguments, or by chasing the diagram around algebraically, we see that some multiple of $a\in H_n(V)$ is in the image of $H_n(F_0)$. So $a$ goes to a torsion element in $\cok(i)$ and so $0$ in $\cok(i)$ mod torsion.

\paragraph{Case: $\mathbf{H_n(V)\ncong H_n(V')}$.} We next consider the case in which $H_n(V)\ncong H_n(V')$. This happens if $\Z\cong H_n(V', V)\to H_{n-1}(V)$ has non-trivial kernel, i.e. if a multiple of the attaching sphere bounds in $V$. In this case, $H_n(V')\cong H_n(V)\oplus \Z$, the additional $\Z$ summand can be taken as generated by a chain $C$ consisting of a multiple of the core of the attached disk $D^n$ and a chain in $V$ whose boundary is a multiple of the attaching sphere. By pushing in along a collar of $\bd V$, we can assume that the geometric intersection of this  chain $C$ with $F$ is the attaching sphere. $C$ is well-defined in this way up to a cycle in $V$, but we can fix a specific one as a generator of the summand.

\subparagraph{Subcase: $\mathbf{H_n(F')\protect\cong H_n(F_0)\oplus \Z}$.}
Suppose that the translate of the attaching sphere, $\bd( D^n\times *)$, $*\in S^{n-1}$, weakly bounds in $F_0$. Then from the long exact sequence of the pair $(F',F_0)$, we see that $H_n(F')\cong H_n(F_0)\oplus \Z$. This follows since $H_n(F',F_0)\cong H_n(D^n\times S^{n-1}, S^{n-1}\times S^{n-1})\cong \Z$, using excision and the long exact sequence of the latter pair. The distinguished $\Z$ summand of $H_n(F')\cong H_n(F)\oplus \Z$ can then be generated by a chain $B$ composed of a multiple of a translate of the core of the handle and another chain in $F_0$ whose boundary coincides with that of this multiple of the core. $B$ is well-defined up to cycles in $F_0$,  and again we fix a representative. The image in $H_n(V',V)$ of the chain $B$ represents a non-trivial multiple of the generator. 

We will study  $\cok(i')$ and $\cok(i_0)$ modulo torsion. Writing $H_n(F')\cong H_n(F_0)\oplus \Z$ and $H_n(V')\cong H_n(V)\oplus \Z$, we have clearly that $i'(x,0)=(i_0(x), 0)$, since the image of $F_0$ is in $V$ and hence all such elements go to $0$ under the surjection $H_n(V')\to \Z\cong H_n(V',V)$. We also have that $i'(0,B)=(y, z)$, where $y$ is unknown at this point, but $z$ must be non-zero, since, again, we know that $B$ represents a non-trivial multiple of the generator of  $H_n(V',V)$. 

Consider the diagram
\begin{diagram}
0&\rTo& H_n(F_0) &\rTo& H_n(F') &\rTo& \Z\cong H_n(F')/H_n(F_0) &\rTo& 0\\
&&\dTo^{i_0} && \dTo^{i'} && \dTo\\
0&\rTo& H_n(V) &\rTo& H_n(V') &\rTo& \Z\cong H_n(V')/H_n(V) &\rTo& 0&.
\end{diagram}
As noted, the righthand vertical map must be injective. Thus its kernel is $0$, and by the serpent lemma the map $\cok(i_0)\to \cok(i')$ is an injection. 

If we consider this diagram with $\Q$ coefficients, the righthand map is also surjective and $\cok_{\Q}(i_0)\cong \cok_{\Q}(i')$, induced by inclusion. 

Now let's look at $H_n(F_0;\Q)\to H_n(F;\Q)$. This is also an injection by the long exact sequence of the pair. Suppose it is not an isomorphism. Then from the long exact sequence of the pair, $H_n(F;\Q)\cong H_n(F_0;\Q)\oplus \Q$. A generator $A$ of the distinguished $\Q$ can be represented by a chain contained in $F$ consisting of a multiple of a fiber of the tubular neighborhood of the attaching disk plus a chain in $F_0$ with the opposite boundary. This is because the existence of this extra term implies that a multiple of the boundary of the fiber bounds in $F_0$. We will see that this situation actually can't arise.

In $H_n(V',F')$, the image of $A$ is clearly 
homologous to a multiple of the relative cycle generated by the cocore of 
the handle, and, by the assumptions of this case leading to the 
non-triviality and non-torsion of $C$, the intersection of $A$ and $C$ 
cannot be $0$, and it would follow that this image of $A$ generates an 
infinite cyclic subgroup of $H_n(V',F';\Q)$. So under the maps $H_n(F)\to H_n(V)\to H_n(V')\to H_n(V',F')$, $A$ must map to a non-trivial element. 
Thus $A$ maps to some element $0\neq x\in H_n(V;\Q)$, which maps to $0\neq (x,0)\in H_n(V';\Q)$. Now consider the image  of $x$ in $H_n(V',F')$. This elements is still represented by $A$, modulo chains in $F'$. The intersection of $A$ with any cycle in $V$ is $0$, since any such cycle can be pushed into the interior of $V$ and thus be made disjoint from $F$ and $F'$. Now consider the intersection of $A$ with $C$. We know that $i'(0,B)=(y,z)$, where $z=mC$ for some $m\in \Q$. But then the intersection of $A$ with $(y,z)$ is $0$, since $(y,z)$ goes to $0$ in $H_n(V',F')$ and since $A$ is the image of an element of $H_n(V')$. But this implies that the intersection of $A$ with $y$ is the negative of its intersection with $z$. But the intersection of $A$ with $y$ is $0$ since $y$ is in $H_n(V)$. Thus the intersection of $A$ with $z$ is $0$, and so the intersection of $A$ with $C$ is $0$. It then follows that $A$ must map to $0$ in $H_n(V',F';\Q)$ since 
$H_n(V',F';\Q)$ and $H_n(V';\Q)$ are dual under the intersection pairing. 
So we arrive at a contradiction. 
Thus it must be in fact that $H_n(F)\cong H_n(F_0)$. 

So we see that that $\cok(i)\cong \cok(i_0)$. However, we still have that 
$\cok (i_0)\to \cok(i')$ may only be an injection, the cokernel of \emph{this} 
map being a cyclic torsion group. We can assume by changing 
basis if necessary that, modulo torsion, this map is represented by a 
matrix that is $0$ except on the diagonal, all diagonal entries except 
perhaps the last one being equal to $1$. The last entry is non-zero, say $p$, but 
may not be $1$. So now all other basis elements of $\cok(i')$ but the last are represented by the chains 
that represent them in $\cok (i)$ mod torsion and so their linking 
pairings with each other remain unchanged. The last basis element is 
homologous to $1/p$ times a chain lying in $\cok (i)$. So each of its 
linking numbers will simply be $1/p$ times those for the corresponding chain 
in $\cok (i)$. Hence the change to the Seifert matrix is to multiply the 
last row and column by $1/p$. In other words, the Seifert matrix changes by a rational change of bases, although the new matrix must also be integral.

\subparagraph{Subcase: $\mathbf{H_n(F_0)\protect\cong H_n(F') \protect\cong H_n(F)}$.}

Suppose $H_n(F_0)\cong H_n(F')$. In this case, we show first that it is impossible to also have $H_n(F_0)\cong H_n(F)$, induced by inclusion. So suppose that $H_n(F_0)\cong H_n(F')\cong H_n(F)$, both isomorphisms induced by inclusion of $F_0$. Then the attaching sphere must generate a torsion (or zero) subgroup of $H_{n-1}(F)$. This is because all cycles of $H_n(F)$ can be homotoped into the interior of $F_0$ so that the intersection of the attaching sphere with any such cycle is empty. Thus, by the Poincar\'e duality of the $2n-1$ manifold $\bd V$, whose homology in all but the top dimension is equal to the homology of $F$, the attaching sphere cannot generate a free subgroup of $H_{n-1}(F)$. It follows that some multiple of the attaching sphere must bound in $F$. Thus, in rational homology, in which $H_{n}(V';\Q)\cong H_n(V;\Q)\oplus \Q$, the distinguished $\Q$ summand can be taken as generated by a cycle $C$ composed of the attaching disk and a chain in $F$ whose boundary is the (negative of) the attaching sphere. A multiple of $C$ will generate the corresponding distinguished $\Z$ term with $\Z$ coefficients. 

Okay, so now if $H_n(F_0)\cong H_n(F)$, $\cok(i)\cong \cok(i_0)$, integrally or rationally and generated by the same cycles in $F_0$. And since $H_n(F_0)\cong H_n(F')$, also generated by the same cycles, $\im(i_0)= \im(i')\subset H_n(V)\subset H_n(V')$, so we see that $\cok_{\Q}(i')\cong \cok_{\Q}(i_0)\oplus \Q$, the distinguished $\Q$ summand again generated by $C$. So the rational Seifert matrix for $V'$ has one more row and column than that for $V$, and except for this row and column is identical to that for $V$. In this row and column, all except possibly the diagonal entry must be $0$ because $C$ cannot link any element in $V$. This is because in the process of putting a cobordism on $F$, we have extended the knot originally in $D^{2n+1}$ to be in $D^{2n+1}\cup S^{2n}\times I$. The cobordism from $F$ lies in $S^{2n}\times I$, and hence so does $C$. But all element representing cycles from $H_n(V)$ lie in the original $D^{2n+1}$. Since the $n$-dimensional homology groups of both $D^{2n+1}$ and $S^{2n}\times I$ are trivial, cycles in each can bound entirely within each (and we can push along some collars if necessary). So $C$ need not link anything from $H_n(V)$. Thus the rational Seifert matrix is $0$ along the additional row and column except where they meet. 

But now this must violate the invariance of the Alexander polynomial, which can be computed from the rational Seifert matrix. If the diagonal term is $0$ or if $n$ is odd, then $R= -\theta'+(-1)^{n+1}\theta$ is singular, which is impossible. If the diagonal term is not $0$, say it is $x\neq 0$, then the Alexander polynomial will be altered by multiplication  by $\frac{xt+ x}{2x}=\frac{t+1}{2}$, which is also impossible as this term is not a rational multiple of a power of $t$ and hence not a unit in the ring of rational Laurent polynomials.

\subparagraph{ Subcase: $\mathbf{H_n(F')\protect\cong H_n(F_0)} $ but  $\mathbf{H_n(F_0)\protect\ncong H_n(F) }$.} 

In this case, $H_n(F)\cong H_n(F_0)\oplus \Z$, from the long exact sequence of $(F,F_0)$. The $\Z$ term can be taken as generated by a chain $A$ that is the sum of a multiple of the fiber disk of the tubular neighborhood of the attaching sphere and another chain in $F_0$ with the opposite boundary. 

The chain $A$ must generate an infinite cyclic summand in $H_n(V)$ because, under the composition $H_n(F)\to H_n(V)\to H_n(V')\to H_n(V',F')$, $A$ becomes relatively homologous to a multiple of the cocore of the attached handle, and this cocore must have a non-zero intersection number with any chain generating the distinguished $\Z$ summand of $H_n(V')\cong H_n(V)\oplus \Z$. We do not here run into the contradiction of the previous similar case since it is no longer true that a multiple of the generator of this summand of $H_n(V')$ is in the image of  $i'$, since now the image of $i'$ in $H_n(V')$  must equal the image of $i_0$ in $H_n(V)\subset H_n(V')$.  Meanwhile, the image of $A$ in $H_n(V)$ must not be in the image of $H_n(F_0)$, since the composition $H_n(F_0)\cong H_n(F')\to H_n(V')\to H_n(V',F')$ is $0$, and we know that the image of $H_n(F_0)$ in $H_n(V')$ is the same as the image of $H_n(F_0)$ in $H_n(V)\subset H_n(V')$. So we see that in fact $A$ generate an infinite cyclic group in $H_n(V)$ that is not in the image of $H_n(F_0)$. So, mod torsion, $\cok(i_0)\cong \cok(i)\oplus \Z$. 

It also follows from the serpent lemma that $\cok(i')\cong \cok (i_0)\oplus \Z\cong \cok(i)\oplus \Z^2$. 

\begin{diagram}
0&\rTo& H_n(F_0) &\rTo& H_n(F') &\rTo& 0&\rTo&0\\
&&\dTo^{i_0} &&\dTo^{i}& &\dTo&\\ 
0&\rTo& H_n(V) &\rTo& H_n(V') &\rTo& \Z&\rTo&0
\end{diagram}

Thus we see that the Seifert matrix for $V'$ has two more rows and columns than the one for $V$, and, excluding these rows and columns, the matrices agree. We must now determine what entries go in these last two rows and columns for $V'$. By changing bases if necessary, we can assume that $A$ is a multiple of a generator of the distinguished $\Z$ term of $\cok(i_0)\cong \cok(i)\oplus \Z$. But as in the previous case, we see that $A$, because it lies in $F$, does not link with any of the cycles in $H_n(V)$ including itself. It can only possibly link nontrivially with a chain generating the distinguished $\Z$ summand of $H_n(V')\cong H_n(V)\oplus \Z$. The same is then true for the generator of the summand containing $A$. Thus the matrix for $V'$ must differ from that for $V$ as in equation \eqref{E: S-equiv}. The same arguments then show that we must have a rational S-equivalence.

\subsubsection[Handles of index $n+1$]{Handles of index $\mathbf{n+1}$}

Consider again the long exact sequence for $(F,F_0)$. By excision, $H_i(F,F_0)\cong H_i(S^n\times D^{n-1}, S^n\times S^{n-2})$. Clearly, $H_{n+1}(S^n\times D^{n-1})=H_{n-1}(S^n\times S^{n-2})=0$, and furthermore, $H_n(S^n\times S^{n-2})\cong H_n(S^n\times D^{n-1})\cong \Z$, the isomorphism being induced by inclusion and taking a generator $S^n\times *\subset S^n\times S^{n-2}$ to a generator $S^n\times *\subset S^n\times D^{n-1}$. It follows that $H_i(S^n\times D^{n-1}, S^n\times S^{n-2})$ and hence $H_i(F,F_0)$ is $0$ for $i=n,n+1$. Thus $H_n(F_0)\cong H_n(F)$, induced by inclusion. Thus from the commutative diagram

\begin{equation*}
\begin{CD}
H_n(F)&@>i>>& H_n(V)\\
@A\cong AA & & @AA=A\\
H_n(F_0)&@>i_0>>& H_n(V),
\end{CD}
\end{equation*}
we see that $\cok(i)=\cok(i_0)$.

On the other hand, we consider the Mayer-Vietoris sequence for $F'$ and $F_0$. Since $H_{n-1}(S^{n}\times S^{n-2})=H_n(D^{n+1}\times S^{n-2})=0$, the inclusion-induced homomorphism $H_n(F_0)\to H_n(F)$ is onto, possibly with kernel represented by the attaching sphere, appropriately translated to $S^n\times *\subset S^n\times S^{n-2}\subset F_0$. 

Meanwhile, since $V'$ is obtained from $V$ by attaching an $n+1$ handle, $H_i(V,V')$ is $0$ for $i\neq n+1$ and $\Z$ for $i=n+1$. Thus $H_n(V)\to H_n(V')$ is also onto, and its kernel is also generated by the attaching sphere. If the class of the attaching sphere is either trivial or torsion in $H_n(V)$, then $H_n(V)\to H_n(V')$ is an isomorphism mod torsion, and we obtain a diagram

\begin{equation*}
\begin{CD}
H_n(F_0)&@>i_0>>& F_n(V)\\
@V\text{onto} VV & & @VV\cong V\\
H_n(F')&@>i'>>& F_n(V').
\end{CD}
\end{equation*}
Again we see that $\cok(i_0)\cong \cok (i)$, and again, since all maps are induced by inclusions, the Seifert pairing is unchanged. 

If the attaching sphere generates an infinite cyclic subgroup of $H_n(V)$, it must also generate an infinite cyclic subgroup of $H_n(F_0)$ (if some multiple of it bounds in $F_0$, then that multiple also bounds in $V$ since $F_0\subset V$). So we have the following diagram
\begin{diagram}
0&\rTo &\Z &\rTo& H_n(F_0) &\rTo& H_n(F') &\rTo &0\\
&&\dTo^{\cong} &&\dTo^{i_0} &&\dTo^{i'}&& \\
0&\rTo &\Z &\rTo & H_n(V) &\rTo & H_n(V') &\rTo &0,
\end{diagram}
in which both $\Z$ summands are generated by the attaching sphere. It follows now from the serpent lemma that $\cok (i_0)\cong \cok (i')$. It once more follows that the Seifert matrix is unchanged.

\bibliographystyle{amsplain}
\bibliography{bib}

Several diagrams in this paper were typeset using the\TeX\, commutative
diagrams package by Paul Taylor.

\end{document}